# Multi-parametric Analysis for Mixed Integer Linear Programming: An Application to Transmission Planning and Congestion Control

Jian Liu[1], Rui Bo[1*], Siyuan Wang[2,1]

[1]Department of Electrical and Computer Engineering, Missouri University of Science and Technology, Rolla, MO 65409, USA

[2]Whiting School of Engineering, Johns Hopkins University, Baltimore, MD 21218, USA

{[1]jliu@mst.edu, [1*]rbo@mst.edu, [2]siyuanwang@jhu.edu}

**Abstract:** Enhancing existing transmission lines is a useful tool to combat transmission congestion and guarantee transmission security with increasing demand and boosting the renewable energy source. This study concerns the selection of lines whose capacity should be expanded and by how much from the perspective of independent system operator (ISO) to minimize the system cost with the consideration of transmission line constraints and electricity generation and demand balance conditions, and incorporating ramp-up and startup ramp rates, shutdown ramp rates, ramp-down rate limits and minimum up and minimum down times. For that purpose, we develop the ISO unit commitment and economic dispatch model and show it as a right-hand side uncertainty multiple parametric analysis for the mixed integer linear programming (MILP) problem. We first relax the binary variable to continuous variables and employ the Lagrange method and Karush-Kuhn-Tucker conditions to obtain optimal solutions (optimal decision variables and objective function) and critical regions associated with active and inactive constraints. Further, we extend the traditional branch and bound method for the large-scale MILP problem by determining the upper bound of the problem at each node, then comparing the difference between the upper and lower bounds and reaching the approximate optimal solution within the decision makers' tolerated error range. In additional, the objective function's first derivative on the parameters of each line is used to inform the selection of lines to ease congestion and maximize social welfare. Finally, the amount of capacity upgrade will be chosen by balancing the cost-reduction rate of the objective function on parameters and the cost of the line upgrade. Our findings are supported by numerical simulation and provide transmission line planners with decision-making guidance.

## 1. Introduction

Many governments and scientists prefer renewable energy to traditional fossil fuels because it produces clean electricity while producing few greenhouse gases or pollutants. Currently, renewable energy accounts for roughly one-quarter of our global electrical supply (EIA (U.S. Energy Information Administration), 2021). Hydropower is by far the most common, accounting for 83% of global renewable energy generation because the technology to generate electricity from water flow has existed for a long time. Wind energy came in second, accounting for little more than 7% of renewable electricity output, followed by biowaste and biomass energy (7%), geothermal energy (2%), and solar, tidal, and wave energy (1%).

In addition to the ongoing increase in demand for renewable energy, renewable energy power generation bears the characteristics of uncertainty (Liu and Sun et al., 2022; Liu and Ou et al., 2022). Renewable energy, poses new challenges, not the least of which is high volatility, which enhances supply-side uncertainty. Renewable energy sources such as solar and wind are often distributed along the grid, producing variable and intermittent power (Bird et al., 2013). The inability of the grid to accommodate a power source with such characteristics can lead to many technical problems such as current congestion, overvoltage, and stability issues (Brown et al., 2012). Furthermore, with increased demand for power, present lines' capacity is straining to keep up, resulting in an evident problem: line congestion. Uncertainty in renewable energy generation and a lack of network capacity may cause a sudden increase in costs and jeopardizes the security of electricity availability.

To handle network congestion while maintaining system security and stability, demand response, transmission expansion planning (TEP), and modifying the transmission capacity of existing lines can be employed. The TEP problem is a study of identifying new transmission lines that need to be added to an existing electrical system. Aiming at the uncertainty of wind turbine load and output power, Ugranli and Karatepe (2016) proposed a new transmission expansion planning methodology considering N-1 contingency conditions and used a cluster-based probabilistic method to determine load and wind models. Muñoz et al. (2019) developed a new method for computing models based on dynamic thermal levels to maximize line current capacity using an iterative span recognition procedure based on restricted conductor deviation. González et al. (2022) extended a multi-level generation and TEP optimization approach to analyze the impact of the emergence of the community energy sector on power system expansion. The cost of increasing the power grid could drastically shift, especially if massive new transmission capacities support renewable energy integration (Sebastian et al., 2022).

Given the many changes in the way the power transmission system is being planned and operated, Larruskain et al. (2006) proposed reaching higher current densities in existing transmission lines to upgrade transmission capacity. Increasing power carrying capability and loss reduction between the conversion of AC distribution lines and DC lines can upgrade transmission capacity. (Larruskain et al., 2011; Manickam, 2018). It is estimated that by 2026, the global investment in global transmission and distribution infrastructure is expected to reach US $ 351 billion, and some of them will eliminate line capacity restrictions (Mbuli et al., 2019). The expansion of the transmission capacity through the use of existing lines and the right to pass is a strategy worth paying more attention to. The demand for the location approval and line of the new line may increase the project schedule for many years, and in some cases, it may even lead to the cancellation of the project. Liza et al. (2020) have developed a classification method of technical strategy, which can be used to increase the capacity of existing transmission corridors and determine the supervision and other issues that must be solved in each situation. As a result, increasing current line transmission capacity is critical to meeting the increasingly rising and uncertain renewable energy transmission demands at the lowest possible cost within a limited budget.

There are various ways above to alleviate line congestion caused by increased renewable energy demand and uncertainty from line planning expansion and capacity increase, but there is no realistic method to estimate cost with increased line transmission capacity when the entire grid is operating. As a result, we intend to develop a model that is more suitable for actual power grid operators to measure the impact of increased line transmission capacity, obtain accurate, explicit solutions using multi-parametric programming methods, and analyze the relationship between increased capacity (i.e., required budget) and decreased system cost to choose the right line(s) to relieve the pressure of line congestion.

Multi-parametric programming can effectively handle our problem. The parameter domain is set to an adjustable range of increasing line transmission capacity. The parameter analysis can produce an explicit piecewise affine formula for the cost of increasing line capacity and decreasing system costs. At the same time, the changing range of line capacity on different lines can be represented by different parameters, and we can know the cost-related parameter changes to select the most suitable line for capacity adjustment by deriving the system cost explicit expressions obtained for each line concerning parameter changes. As an independent system operator (ISO) in the electricity market, ISO clears energy and operates reserves, commits generators with the lowest price offers based on the required demand to minimize generation costs or maximize social welfare, and then generates the cleared locational marginal pricing (LMP) (Chen and Li, 2011; Soofi and Manshadi, 2022). In contrast to the previous power grid linear economic dispatch programming model (Li, 2007; Li and Bo, 2009), which focused on demand uncertainty, we include a large scale of binary variables in the unit commitment and economic dispatch problem to better reflect transmission line capacity upgrade and ISO's goal-minimizing system cost or maximizing social welfare.

Unlike the prior problem, which involved solving a single period and generating wind speed and load value based on their distribution to calculate Probability Distribution Function for each variable (Hemmati et al., 2014; Ugranli and Karatepe, 2016). Rivotti and Pistikopoulos (2014) studied constrained dynamic programming for multi-stage mixed-integer linear formulations with a linear objective function. The current problem is eventually characterized by a multi-period and multi-parametric mixed-integer linear programming model. The B&B method (Acevedo and Pistikopoulos, 1999) and iterating between two subproblems decomposed by the multi-parametric mixed integer linear programming (mp-MILP) problem (Dua et al., 2001) are traditional methods for multi-parametric mixed integer programming problems. However, the B&B method faces the challenge of effectively solving the problem when there are too many nodes/binary variables to explore because we should check all possible nodes to determine the best solution, which increases to the computation load because the current study only provides the lower boundary in each iteration by relaxing some integer decision variables to continuous variables to perform parametric analysis.

To solve the above challenges, we aim to address the following questions: 1) How can approximately affine solutions to a multi-period and multi-parametric mixed integer linear problem with transmission capacity uncertainty be found efficiently? 2) Given the limited budget, which lines' capacity should be upgraded and how much to enhance societal welfare or reduce system costs? 3) What is the functional relationship between system cost and line transmission capacity adjustment? To that purpose, we develop the right-hand side (RHS) parametric analysis MILP problem model to reflect the uncertainty in line transmission capacity. To begin, solve the relaxed unit commitment and economic dispatch (UCED) problem at each node using the classical branch and bound method and obtain the corresponding lower boundary. Then, using our proposed approach, we will transform the continuous variables into binary variables and fix them so that the upper boundary for each node is achieved. Third, to find the best approximation solution, compare the lower and higher boundaries of each node sequentially; if the error is smaller than the tolerance range, the branching process is ended, and the current upper bound is our final approximate solution. To the best of our knowledge, this is the first work to investigate the upper and lower bounds to achieve the best approximate solution for multiple periods of UCED parametric analysis problems with large binary decision variables.

Our research makes three significant contributions: To begin, we formulate the problem in matrix format and use the language technique and the KKT condition to get the best solution at each node by relax the binary variables (i.e., UC) to continuous variables. We use parametric analysis for MILP to develop a mathematical function link between system cost and expandable line capacity. We present a method for determining the best objective function and related critical region by utilizing active and inactive constraints determined through optimal language multiples.

Second, multiparameter programming issues using integers have exponential complexity, which means that as the number of binary variables increases, so does the computing burden and simulation time. Our approach not only generates the exact optimal analytical solution for the small-scale mp-MILP parametric analysis problem but also produces the best approximation analytical results for the large-scale mp-MILP parametric analysis problem by calculating the upper and lower boundaries at each node. This is a valuable extension to the traditional B&B algorithm.

Third, enhancing existing transmission line capacity is a beneficial strategy to respond to development by offsetting instability and maintaining the security of transmission lines with increased demand and the development of high uncertainty renewable energy sources. We investigate the functional relationship between system cost and increased transmission capacity for multiple periods of ISO's UCED problem utilizing RHS uncertainty of parametric analysis. Analyzing the first order function between the objective function and the parameters enables decision makers to implement suitable line capacity expansion and achieve the best social welfare with the least amount of investment.

The remainder of the work is structured as follows. Section 2 examines relevant studies. Section 3 establishes the model and formulates the transmission line capacity expansion problem from the perspective of ISO utilizing the UCED problem. Section 4 discusses RHS parametric analysis for large-scale MILP optimization. Section 5 evaluates the proposed approach by applying it to a single-period problem and a 24-period IEEE-5 bus UCED problem. Finally, section 6 summarizes the work and provides recommendations for further research.

## 2. Literature review

This paper discusses the implications of multi-parametric programming models, approaches for grappling with multi-parametric mixed integer programming problems, and the electricity market with uncertainty. Each of these is reviewed separately.

### *2.1 Multi-parametric Programming Model*

Multi-parametric programming implications are becoming widely used in operations research. Dua et al. (2001) introduced a model that considers the control variables as optimization variables and the state variables as parameters to handle engineering problems reformulated as multi-parametric mixed-integer quadratic problems (mp-MIQPs) in mixed logical dynamical systems. Later, Dua et al. (2002) used the lagrange method to solve multi-parametric quadratic problems (mp-QPs) and mp-MIQPs with a convex and quadratic objective function and linear constraints, using an affine expression for the optimal solution to systematically characterize the space of parameters by a set of optimality regions in model predictive and hybrid control problems. Gupta (2011) suggested a new method for multi-parametric programming problems based on the enumeration of active sets, it avoids the enumeration procedure's combinatorial explosion and makes the enumeration implicit to solve mp-QPs.

Habibi et al. (2016) proposed a multi-parametric programming model for obtaining the explicit solution to optimal control problems for some classes of hybrid systems, as well as an approximation algorithm for solving a general type of mp-MILP problems. Avraamidou et al. (2017) defined seasonal demand variability in supply chain planning problems as multi-parametric mixed-integer bi-level linear programming problems and proposed a novel algorithm that provides an exact, global, and parametric solution with or without demand uncertainty.

Mate et al. (2020) proposed an offline combinatorial approach to identify all active sets of constraints for the nonlinear model predictive control (MPC) problem a priori applying KKT conditions. Mate (2020) investigated KKT conditions to identify active constraints and presented piecewise affine control laws and their accompanying critical regions as part of the multi-parametric MPC technique for linear systems. Shokry et al. (2021) proposed multiparametric programming to solve chemical process operation optimization problems with unavoidable uncertain parameters. Pappas et al. (2021) proposed an algorithm for the exact solution of explicit nonlinear MPC problems with convex quadratic constraints based on a second-order Taylor approximation of Fiacco's Basic Sensitivity Theorem applied to an exact nonlinear MPC problem.

According to the above, it can be known that the multi-parametric programming model can solve the uncertainty problems in many fields. Unlike the previous study, we propose an mp-MILP model in the multi-period electricity market to address the challenges incorporating the uncertainty of transmission line capacity in the UCED problem in the electricity market, which contains binary variables reflecting the unit commitment of generators; we characterize the problem mathematically and obtain an analytical solution.

### *2.2 Multi-parametric Mixed Integer Programming Problems*

Acevedo and Pistikopoulos (1997) presented a novel branch and bound algorithm based on the solution of multi-parametric linear programs at each node of the tree search and special bonding procedures to solve MILP problems where the right-hand-side parameters are allowed to vary respectively to analyze linear process engineering problems under uncertainty. Dua and Pistikopoulos (2000) decomposed the mp-MILP problem into two subproblems and then iterated between them. The first subproblem is formulated as an mp-LP problem by fixing integer variables, and the second subproblem is formulated as a MILP problem by relaxing the parameters as variables. Faísca et.al (2009) developed a method for solving the mp-MILP problem by changing parameters in the objective function and the RHS of the constraints and splitting the mp-MILP problem into two sub-problems, a master MINLP problem and a slave multi-parametric global optimization problem.

In the process of solving mp-MILP, mp-LP problem must be solved firstly. Gal and Nedoma (1972) was the pioneer who considered the computational aspects of the mp-LP problem, presented fundamental notations, theorems, and definitions, and showed how to solve the mp-LP issue. Mitsos and Barton (2009) introduced an enhancement of the well-known rational simplex approach for parametric LPs that involves fewer successive operations on rational functions. Support set invariancy and optimal partition invariancy were presented by Hladík (2010), who compared them to the classical optimal basis approach to solve multi-parametric linear programming problem. Adelgren and Wiecek (2016) suggested a new two-phase method for addressing the multi-parametric linear complementarity problem with sufficient matrices, where the mp-LCP solution is invariant over each region as a function of the parameters.

Ting et al. (2016) used MILP relaxations from multi-parametric disaggregation to construct rigorous relaxations of increasing complexity for the two different nonlinear formulations. Al-Shihabi et al. (2022) proposed a basic variable neighborhood search approach to solve the repatriation scheduling problem using MILP. Nobil et al. (2022) proposed a MILP model with an objective cost function to determine homogeneous personnel shift scheduling so that the number of personnel per shift constraints is met and to optimize personnel scheduling with the same skills to minimize variable hospital costs.

Crema (1997, 2002) created a method for achieving a comprehensive multi-parametric analysis by solving a family of ILP problems and selecting a suitable finite sequence of nonparametric MILP questions. Domínguez and Pistikopoulos (2010) proposed two multiparametric programming algorithms for handling pure integer and mixed-integer bilevel programming problems; the second algorithm targets the mixed-integer case of the bilevel programming problem. Oberdieck et al. (2014) proposed a new algorithm for solving mp-MILP problems that employ a branch and bound strategy and McCormick relaxation procedures to overcome the presence of bilinear terms in the model, resulting in an envelope of parametric profiles containing the optimal solution to the mp-MILP problem.

We specified and spliced constraints in the model with matrices, both equation and inequation constraints were presented in detail. By considering the uncertainty of transmission line capacity, we investigated an effective approach for the parametric analysis in multiple periods of MILP of ISO's UCED scheduling problem, which has not been addressed previously.

### *2.3 Power Market under Uncertainty*

In the United States, both day-ahead and real-time electricity markets have the market-clearing process modeled as a UCED problem (ERCOT, 2020; PJM, 2021). Carrión and Arroyo (2006) presented a new mixed-integer linear formulation requiring fewer binary variables and constraints for the unit commitment problem of thermal units. Lin et al. (2018) described a transmission and distribution network coordinated dynamic economic dispatch (DED) model and proposed an efficient decentralized method to solve this problem using multi-parametric quadratic programming. Wang et al. (2022) proposed a two-stage distributionally robust unit commitment framework with both regular and flexible generation resources.

Resources of a wind-thermal DED problem are fluctuating, depending on the weather conditions (Zaman et.al, 2016). Moarefdoost et al. (2016) provided an alternative way to considering the uncertainty of renewable energy sources and the consequent ramping of conventional generation by a robust reformulation of the problem. Wind and solar are highly intermittent in availability and result in uncertainty in demand fulfillment. The generation of these sources may be different according to the season and weather conditions of the day. For example, changes in the degree of heat can lead to fluctuations in the power output of the solar illumination power; similarly, the target drawn from the wind power system may suddenly change with the wind speed (Pravin et al., 2020). Prajapati et al. (2021) suggested the installation of an energy storage system to deal with the uncertainties causing by the high penetration of renewable energy sources making the power system unreliable.

In the deregulated electricity market, Hemmati et al. (2014) established a novel method for TEP, which is associated with reactive power planning, reliability evaluation, and the consideration of wind and load uncertainty. Because stochastic programming models (SPM) couldn't solve the huge number of scenarios required to appropriately describe the stochastic and non-dispatchable characteristics of uncertain renewable sources, Zugno and Conejo (2015) substituted SPM with an adaptive resilient optimization problem. Ishizaki et al. (2020) created an energy market model based on a robust convex program, as well as modeling and analysis of day-ahead spatial temporal energy markets. Yang et al. (2021) proposed a novel model-free dynamic dispatch strategy for integrated energy systems based on improved deep reinforcement learning to address the issue that most existing dynamic dispatch schemes are limited by forecasting or model accuracy due to the randomness of renewable energy generation and demand.

From the perspective of the ISO, modeling and optimizing generators across multiple market clearing processes with uncertainties and incomplete information, such as load demand, raises new challenges. Li (2007) proposed an approach to eliminate the step change in the continuous locational marginal price curve with respect to load fluctuation and smooth the price curve step changes. Li and Bo (2009) introduced a more efficient algorithm to identify the new binding constraint and the new marginal unit set when the load increases. Vaskovskaya et al. (2018) employed a non-linear alternating current power system model to develop analytical formulae for expressing LMPs. Kara et al. (2022) proposed a stochastic local flexibility market to solve grid issues such as voltage deviations and grid congestion in a distribution grid under demand uncertainty and random bidding process.

Unlike the previous work, we first find congested lines by exploiting the relationship between the optimal solution and the uncertain parameters of line capacity, and then use parametric analysis to assess the impact of transmission line capacity increases on optimization targets. We also detect which lines' capacity should be improved and by how much, considerably boosting the social welfare of the system (i.e., maximizing social welfare through minimum investment).

## 3 Model Setup and Analysis

We made some assumptions based on Li (2007), Li and Bo (2009), and Vaskovskaya et al. (2018) to help clarify our model. They are listed as follows:

1) Each bus has one generator and one load for convenience.

2) The cost of electricity generation is linear.

With reference to previous papers, we constructed the following model, which is an mp-MILP model in equation (1).

$$\min \sum_{t=1}^{T}\sum_{i=1}^{N} C_{it}\times G_{it} + UC_{it}\times U_{it} + SC_{it}\times V_{it} \tag{1}$$

$$st.\begin{cases}
\sum_{i=1}^{N} G_{it} - \sum_{i=1}^{N} D_{it} = 0,\ t=1,2,\cdots,T; for\ i=1,2,\cdots,N & (1a)\\
\sum_{i=1}^{N} GSF_{k-i}\times\left(G_{it}-D_{it}\right) \le F_k^{\max}+\theta_k,\ \ t=1,2,\cdots,T; for\ k=1,2,\cdots K & (1b)\\
G_{it}-G_{i(t-1)}-(\mathrm{SR}_i-G_{it}^{\max})\cdot U_{it}-(UR_i-\mathrm{SR}_i)\cdot U_{i(t-1)} \le G_{it}^{\max},\ t=1,2,\cdots,T; for\ i=1,2,\cdots,N & (1c)\\
G_{it}-(G_{it}^{\max}-\mathrm{SD}_i)\cdot U_{i(t+1)}-\mathrm{SD}_i\cdot U_{it}\le 0,\ t=1,2,\cdots,T-1; for\ i=1,2,\cdots,N & (1d)\\
G_{i(t-1)}-G_{it}-(DR_i-SD_i)\cdot U_{it}-(SD_i-G_{it}^{\max})\cdot U_{i(t-1)}\le G_{it}^{\max},\ t=1,2,\cdots,T; for\ i=1,2,\cdots,N & (1e)\\
U_{it}-U_{i(t-1)}\le V_{it},\ t=2,3,\cdots,T & (1f)\\
\sum_{k=t-T_i^{on}+1}^{t} V_{ik}-U_{it}\le 0, t=T_i^{on}+1, T_i^{on}+2,\cdots,T & (1g)\\
\sum_{k=t-T_i^{off}+1}^{t} V_{ik}+U_{i(t-T_i^{off})}\le 1, t=T_i^{off}+1, T_i^{off}+2,\cdots,T & (1h)\\
G_{it}^{\min}\cdot U_{it}\le G_{it}\le G_{it}^{\max}\cdot U_{it},\ t=1,2,\cdots,T; for\ i=1,2,\cdots,N & (1i)\\
U_{it}\in\{0,1\}, the\ unit\ commitment\ for\ generator\ i\ in\ period\ t & (1j)\\
V_{it}\in\{0,1\}, the\ unit\ commitment\ for\ generator\ i\ at\ beginning\ of\ period\ t & (1k)
\end{cases}$$

Where, (1*a*) denotes power balance (i.e., energy supply matches the demand), in each period, the sum of the power generation of all generators should be equal to the sum of demand, (1*b*) presents transmission capacity limit on each line. Decision variables are constrained by startup ramp rates and ramp-up rates in (1*c*), as well as by shutdown ramp rates in (1*d*), (1*e*) presents generation dispatch constrained by shutdown ramp rates and ramp-down rates. The relationship between the state of unit commitment (UC) and the UC start state is seen in (1*f*). Equations (1*g*) -(1*h*) denote minimum up and down time limits, $T_i^{on}$ means minimum on (up) time for the i-th generator, $T_i^{off}$ means minimum off (down) time for the i-th generator. (1*i*) is the constraint for the generation dispatch limit corresponding to the UC for the i-th generator in the t-th period.

In contrast to the previous study on parametric analysis in electricity market problems using multi-parametric linear programming, we focus on the impact of transmission line capacity uncertainty, using a mp-MILP method to construct and execute a more realistic model. However, the existing binary variables (i.e., unit commitment) will make the problem more challenging to solve. As a result, the proposed approach in this work shows the affine relationship between the objective function and transmission line capacity change with uncertainty parameters.

$$\begin{cases} N & number\ of\ buses \\ K & number\ of\ lines \\ C_{it} & generation\ \cos t\ at\ Bus\ i\ in\ period\ t(\$/MWh) \\ G_{it} & generation\ dispatch\ at\ Bus\ i\ in\ period\ t(MWh) \\ UC_{it} & fixed\ \cos t\ if\ turn\ on\ at\ Bus\ i\ in\ period\ t(\$/MWh)\ fixed \\ U_{it} & is\ 1\ (if\ turn\ on\ Bus\ i\ at\ in\ period\ t), else\ is\ 0 \\ SC & start\ \cos t\ if\ turn\ on\ at\ Bus\ i\ in\ period\ t(\$/MWh) \\ V_{it} & is\ 1(if\ turn\ on\ Bus\ i\ at\ the\ beginning\ of\ period\ t)\ , or\ is\ 0 \\ D_{it} & Demand\ at\ Bus\ i\ in\ period\ t(MWh) \\ GSF_{k-i} & generation\ shift\ factor\ to\ line\ k\ from\ bus\ i \\ F_k^{\max} & transmission\ \lim it\ of\ line\ k \\ G_{it}^{\max} & \max imum\ power\ generation\ capacity\ at\ Bus\ i\ in\ period\ t \\ SR_i & startup\ ramp\ \lim it\ at\ bus\ i \\ UR_i & ramp-up\ \lim it\ at\ bus\ i \\ SD_i & shutdown\ \lim it\ at\ bus\ i \\ DR_i & ramp-down\ \lim it\ at\ bus\ i \end{cases}$$

## 4 Optimization and Analysis

To solve the RHS parametric analysis for multiple periods MILP problem, following the previous studies (Acevedo and Pistikopoulos,1997; Oberdieck et al., 2014)), the method based on the B&B method is employed. The proposed algorithm and related procedures are illustrated in the following steps, as indicated in Figure 1.

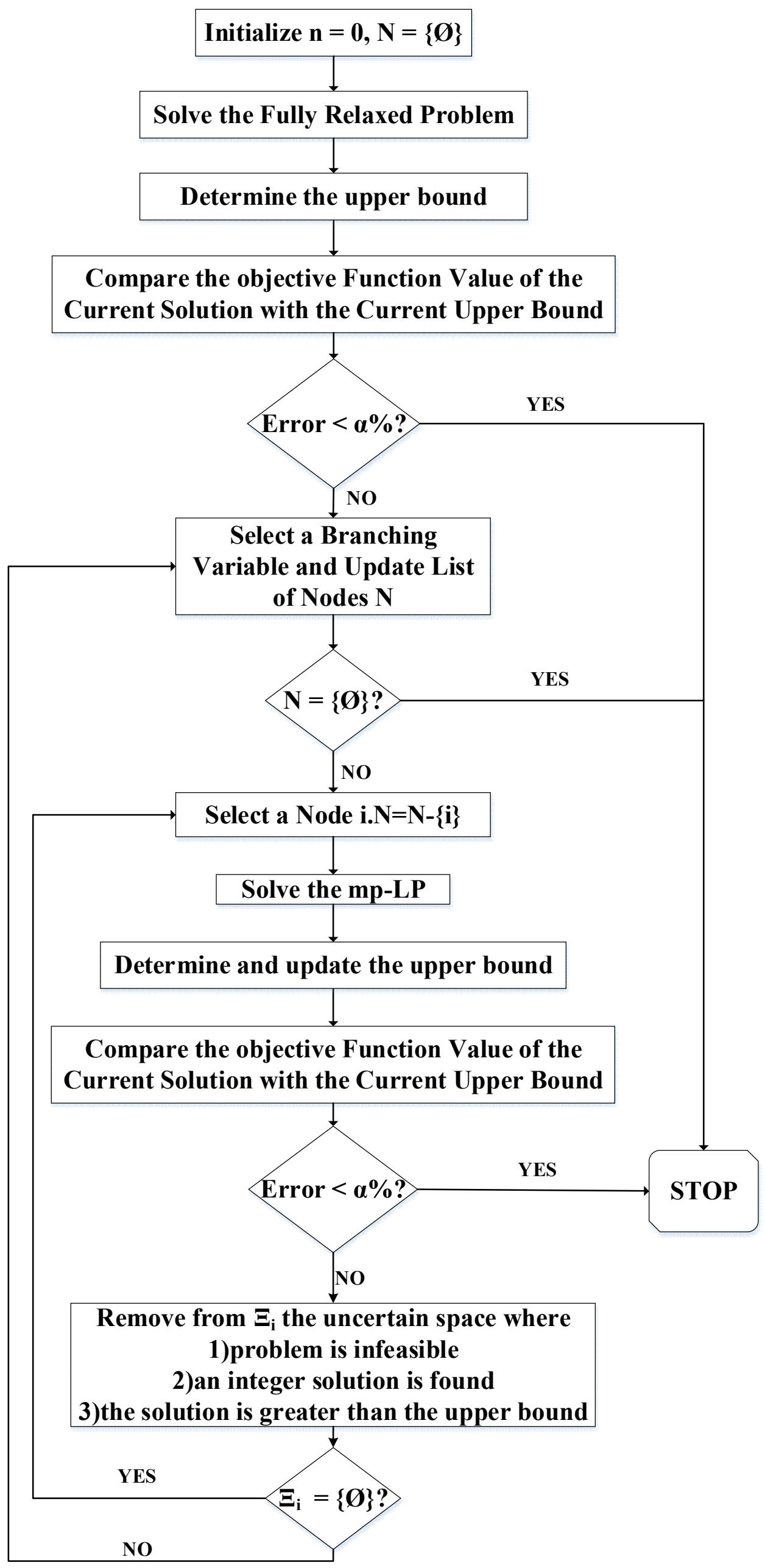

**Figure 1. Algorithm for multi-parametric MILPs**

***Step0. Initialization***

Initialize the collection of nodes that need to be solved, $N = \{\varnothing\}$, recording its counter $n = 0$.

***Step1. Solve the fully relaxed problem***

***1.1 Solve the fully relaxed problem***

After initializing in step zero, a fully relaxed problem must be addressed at the root node to discover the lower bound in the entire parameter space, as shown in Figure 1. The fully relaxed problem widens the range of binary variables. Now, these are continuous variables with values ranging from 0 to 1. As a result, the mp-MILP problem in Eq. (1) is now transforming into a multiple parametric linear programming (mp-LP) problem (2).

$$\min \sum_{t=1}^{T}\sum_{i=1}^{N} C_{it} \times G_{it} + UC_{it} \times U_{it} + SC_{it} \times V_{it}$$

$$st.\begin{cases} Other\ constraints\ are\ the\ same\ with\ (1a) \sim (1i)\ in\ (1) \\ 0 \le U_{it} \le 1 \\ 0 \le V_{it} \le 1 \end{cases} \tag{2}$$

To make the problem in Eq. (2) seem straightforward and facilitate the solution of subsequent problems, in this paper, we write it in a matrix structure as shown in equation (3).

$$z = \min_{G,U,V} M^T \omega$$

$$st.\begin{cases} A\omega \le b + F^{\max} + \theta \\ H\omega = H\hat{D} \end{cases} \tag{3}$$

Following the current work (Oberdieck et al., 2014), we combine the continuous variable representing power generation or economic dispatch and the binary variable representing generator UC, then include both in the same matrix. This will lay a solid foundation for the subsequent solution of the Lagrange Multiplier Method.

Many auxiliary submatrices are defined in Eq. (3), such as the matrix *M*, which is composed of three submatrices made of all the coefficients information of the decision variables in the objective function. Submatrix $C^T$, for example, has NT rows corresponding to the coefficient of per-unit electricity generation cost for all generators in each period, submatrix $UC^T$ has NT rows corresponding to the fixed cost coefficients of UC, and submatrix $\mathrm{SC}^T$ has N(T-1) rows corresponding to the coefficient of the start cost in the objective function.

$$M = \begin{pmatrix} C^T \\ UC^T \\ \mathrm{SC}^T \end{pmatrix}_{N(3T-1)\times 1} = \left(C_{11} \cdots C_{1T} \cdots C_{N1} \cdots C_{NT}\ UC_{11} \cdots UC_{1T} \cdots UC_{N1} \cdots UC_{NT}\ SC_{12} \cdots SC_{1T} \cdots SC_{N2} \cdots SC_{NT}\right)^T \tag{3a}$$

$D, G,$ and $U$ are some fundamental vectors, each with NT rows and 1 column, *D* represents demand for each bus in each period (e.g., $D_{11}$ indicates the demand for bus 1 at period 1), *G* represents generation dispatch decision variables for each bus in each period, and *U* represents unit commitment decision variables for each bus in each period. The startup or shutdown status for each generator is represented by $V$, which is an N(T-1) row and 1 column vector.

$$\begin{aligned} D &= \left(D_{11} \cdots D_{1T} \cdots D_{N1} \cdots D_{NT}\right)^T_{1\times NT} \\ G &= \left(G_{11} \cdots G_{1T} \cdots G_{N1} \cdots G_{NT}\right)^T_{1\times NT} \\ U &= \left(U_{11} \cdots U_{1T} \cdots U_{N1} \cdots U_{NT}\right)^T_{1\times NT} \\ V &= \left(V_{12} \cdots V_{1T} \cdots V_{N2} \cdots V_{NT}\right)^T_{1\times N(T-1)} \end{aligned} \tag{3b}$$

Then, utilizing the aforementioned vectors, next, we will define a new vector $\omega$ in Eq. (3c), which is made up of three sub-vectors of G, U, and V.

$$\omega = \left(G \ U \ V\right)^T_{1\times N(3T-1)} \tag{3c}$$

It is worth noting that A is the left-side coefficients matrix for all of the inequality constants in the Eq. (2), which has $KT + 12NT - 4N - \sum_{i=1}^{N}(T_i^{on} + T_i^{off})$ rows and $N(3T-1)$ columns.

The zero elements in the first row reflect the zero matrices required to keep A intact. Diagonal submatrices are minus one element and one element, respectively. For simplify, in this paper, we use ( $Num$ ) in the subscript to represent $KT + 12NT - 4N - \sum_{i=1}^{N}(T_i^{on} + T_i^{off})$. Appendix A has a detailed summary of the A matrix.

$$A = \begin{pmatrix} GSF & 0 & 0 \\ & \textit{ramp up} & \\ & \textit{shut down} & \\ & \textit{ramp down} & \\ & \textit{state transition} & \\ & \min \textit{on} & \\ & \min \textit{off} & \\ -1 & \hat{G}^{\min} & 0 \\ 1 & -\hat{G}^{\max} & 0 \\ 0 & -1 & 0 \\ 0 & 1 & 0 \\ 0 & 0 & -1 \\ 0 & 0 & 1 \end{pmatrix}_{Num\times N(3T-1)} \tag{3d}$$

In Eq. (3e), $b$ and $\hat{F}^{\max}$ are vectors with $Num$ rows and 1 column. In these vectors, the one and zero elements represent row vector with different sizes.

$$\begin{aligned} b &= \left(GSF\cdot D \ G_i^{\max} \ 0 \ G_i^{\max} \ 0\ 0\ 1\ 0\ 0\ 0\ 1\ 0\ 1\right)^T_{1\times Num} \\ \hat{F}^{\max} &= \left(F^{\max} + \theta \ 0\ 0\ 0\ 0\ 0\ 0\ 0\ 0\ 0\ 0\ 0\ 0\right)^T_{1\times Num} \end{aligned} \tag{3e}$$

In Eq. (3f), $\theta$ and $F^{\max}$ are row vectors with KT rows and 1 column, $\theta$ denotes the uncertainty parameters for each line, which also indicates the capacity that can be enhanced. $F^{\max}$ denotes the current maximum capacity for lines. In this paper, $\theta$ and $F^{\max}$ are two vectors, each with KT rows and 1 column, respectively.

$$\begin{aligned}\theta &= \left(\theta_1 \cdots \theta_1\ \theta_2 \cdots \theta_2 \cdots \theta_K \cdots \theta_K\right)^T_{1\times KT} \\ F^{\max} &= \left(F_1^{\max} \cdots F_1^{\max}\ F_2^{\max} \cdots F_2^{\max} \cdots F_K^{\max} \cdots F_K^{\max}\right)^T_{1\times KT}\end{aligned} \tag{3f}$$

Using Eq. (3), we obtain our matrix structure requirements for the equality constraints from H, which is a matrix of T rows and N(3T-1) columns.

$$H\omega = \begin{pmatrix} \underbrace{1\cdots 0}_{|N|} & \underbrace{1\cdots 0}_{|N|} & \cdots\cdots\cdots & \underbrace{1\cdots 0}_{|N|} & 0 \cdots\cdots\cdots 0 \\ \underbrace{0\ 1\cdots 0}_{|N|} & \underbrace{0\ 1\cdots 0}_{|N|} & \cdots & \underbrace{0\ 1\cdots 0}_{|N|} & 0 \cdots\cdots\cdots 0 \\ \vdots & & \vdots & & \vdots \\ \underbrace{0\cdots 1}_{|N|} & \underbrace{0\cdots 1}_{|N|} & \cdots\cdots\cdots & \underbrace{0\cdots 1}_{|N|} & \underbrace{0 \cdots\cdots\cdots 0}_{N(2T-1)} \end{pmatrix}_{T\times N(3T-1)} \begin{pmatrix} G \\ \bar{U} \\ V \end{pmatrix}_{N(3T-1)\times 1} = \left(\sum_{i=1}^{N} G_{it}\right) \tag{4a}$$

$$H\hat{D} = \begin{pmatrix} \underbrace{1\cdots 0}_{|N|} & \underbrace{1\cdots 0}_{|N|} & \cdots\cdots\cdots & \underbrace{1\cdots 0}_{|N|} & 0 \cdots\cdots\cdots 0 \\ \underbrace{0\ 1\cdots 0}_{|N|} & \underbrace{0\ 1\cdots 0}_{|N|} & \cdots & \underbrace{0\ 1\cdots 0}_{|N|} & 0 \cdots\cdots\cdots 0 \\ \vdots & & \vdots & & \vdots \\ \underbrace{0\cdots 1}_{|N|} & \underbrace{0\cdots 1}_{|N|} & \cdots\cdots\cdots & \underbrace{0\cdots 1}_{|N|} & \underbrace{0 \cdots\cdots\cdots 0}_{N(2T-1)} \end{pmatrix}_{T\times N(3T-1)} \begin{pmatrix} D \\ 0 \end{pmatrix}_{N(3T-1)\times 1} = \left(\sum_{i=1}^{N} D_{it}\right) \tag{4b}$$

Recall the process of Lagrange method and related conditions for linear programming problems with both inequality and equation constraints. Generally, we characterize a linear programming problem as follows:

$$\begin{aligned}&\min f(x) \\ &s.t. \begin{cases} g_i(x) \le 0, for\ all\ i = 1,2,\cdots,m \\ h_j(x) = 0, for\ all\ j = 1,2,\cdots,n \end{cases}\end{aligned} \tag{5}$$

After constructing the Lagrange function of Eq. (5), the KKT conditions are shown as

$$\frac{\partial L}{\partial x} = \frac{\partial f(x^*)}{\partial x} + \sum_{i=1}^{m} \mu_i \cdot \frac{\partial g_i(x^*)}{\partial x} + \sum_{i=1}^{n} \lambda_i \cdot \frac{\partial h_j(x^*)}{\partial x} = 0 \tag{6a}$$

$$g_i(x^*) \le 0, \mu_i \ge 0, for\ all\ i = 1,2,\cdots,m \tag{6b}$$

$$\mu_i g_i(x^*) = 0, for\ all\ i = 1,2,\cdots,m \tag{6c}$$

$$h_j(x^*) = 0, for\ all\ j = 1,2,\cdots,n \tag{6d}$$

In our model, we can rephrase the original problem in Eq. (3) in the following style.

$$\begin{aligned} &\min_{G,U} M^T\omega \\ &st. \begin{cases} A\omega - b - F^{\max} - \theta \le 0 \\ H\omega - H\hat{D} = 0 \end{cases} \end{aligned} \tag{7}$$

Based on Eq. (7), we draw the constraints as follows, assuming $g_i(\omega)$ denotes an inequality constraint and $h_j(\omega)$ stands for an equality constraint.

$$\begin{cases} g_i(\omega) = A_i\omega - b_i - F_i^{\max} \\ h_j(\omega) = H_j\omega - H_j\hat{D} \end{cases} \tag{8}$$

Here, the subscript i and j denote the i-th and j-th row of the corresponding matrices, respectively. We shall take the formulas in Eq. (7) in the following form in accordance with the Lagrange method standards.

$$L = f(\omega) + \sum_{i=1}^{Num} \mu_i \cdot g_i(\omega) + \sum_{j=1}^{T} \lambda_j h_j(\omega) \tag{9}$$

$$f(\omega) = M^T\omega \tag{9a}$$

$$g(\omega) = A\omega - b - \hat{F}^{\max} \le 0 \tag{9b}$$

$$h(\omega) = H\omega - H\hat{D} = 0 \tag{9c}$$

Suppose $\omega^*$ is the optimal solution to the problem, the best response function (i.e., the first-order derivative function) of equation (9) as shown as

$$\begin{aligned} &\frac{\partial L}{\partial \omega} = \frac{\partial f(\omega^*)}{\partial \omega} + \sum_{i=1}^{Num} \mu_i \cdot \frac{\partial g_i(\omega^*)}{\partial \omega} + \sum_{j=1}^{T} \lambda_j \cdot \frac{\partial h_j(\omega^*)}{\partial \omega} \\ &= M + A^T\mu + H^T\lambda = 0 \\ &where \\ &\begin{cases} f(\omega^*) = M^T\omega^*, \nabla f(\omega^*) = M \\ g_i(\omega^*) = A_i\omega^* - b_i - \hat{F}_i^{\max}, \nabla f(\omega^*) = A^T \\ h_j(\omega^*) = H_j\omega^* - H_j\hat{D}, \nabla h(\omega^*) = H^T \end{cases} \end{aligned} \tag{10}$$

We can get the following findings based on the KKT conditions described above.

$$\nabla_x L = \frac{\partial L}{\partial \omega} = 0 \Rightarrow M + A^T\mu + H^T\lambda = 0 \tag{11a}$$

$$\nabla_{\lambda_j} L = \frac{\partial L}{\partial \lambda_j} = 0 \Rightarrow H\omega - H\hat{D} = 0, j = 1,2,\ldots T \tag{11b}$$

$$g_i(\omega) = A\omega - b - \hat{F}^{\max} \le 0, i = 1,\ldots, Num \tag{11c}$$

$$\mu_i \geq 0, i = 1,2,...,Num \tag{11d}$$

$$\mu_i g_i(\omega) = \mu^T \left( A\omega - b - \hat{F}^{\max} \right) = 0, i = 1,...,Num \tag{11e}$$

By Eq.(11b) and (11e), we can get:

$$\mu^T \left( A\omega - b - \hat{F}^{\max} \right) = 0 \tag{12a}$$

$$\lambda^T (H\omega - H\hat{D}) = 0 \tag{12b}$$

Because the convex LP problem's best solutions are always located at the vertices of its feasible region, so, the matrix A will be split into two submatrices, $A_{p_1}$ and $A_{s_1}$, where $p_1$ and $s_1$ represent the index sets associated to the active and inactive constraints. It can be expressed as $p_1 + s_1 = Num$.

$$\mu_{p_1}^T \left( A_{p_1}\omega - b_{p_1} - \hat{F}_{p_1}^{\max} \right) + \mu_{s_1}^T \left( A_{s_1}\omega - b_{s_1} - \hat{F}_{s_1}^{\max} \right) = 0 \tag{13}$$

For inactive constraints in Eq. (13), we also have Lagrange multipliers of $\mu_{s_1} = 0, \mu_{p_1} \neq 0$. When degeneration occurs, the coefficient of the equality constraint may also be 0, so we also divide the equality constraint into two parts. It can be expressed as $p_2 + s_2 = T, \lambda_{s_2} = 0, \lambda_{p_2} \neq 0$

$$\lambda_{p_2}^T (H_{p_2}\omega - H_{p_2}\hat{D}) + \lambda_{s_2}^T (H_{s_2}\omega - H_{s_2}\hat{D}) = 0 \tag{14}$$

Let: $A_p = \begin{pmatrix} A_{p_1} \\ H_{p_2} \end{pmatrix}, b_p + \hat{F}_p^{\max} = \begin{pmatrix} b_{p_1} + \hat{F}_{p_1}^{\max} \\ H_{p_2}\hat{D} \end{pmatrix}$. Then, $A_p$ is a square and invertible matrix; here, $p = N(3 \cdot T - 1)$, so we will get the following results as shown in Eq. (15).

$$A_p\omega^* - b_p - \hat{F}_p^{\max} = 0 \Rightarrow \omega^* = A_p^{-1}(b_p + \hat{F}_p^{\max}) \tag{15}$$

In Eq. (15), let the set $S$ represents the subscripts collection of active constraints from $A$ for subsequent calculations and analysis. Plugging $\omega^*$ into $f(\omega)$ of Eq. (7), we get the following optimal solution.

$$\begin{aligned} z^*(\theta) &= M^T\omega^* \\ &= M^T A_p^{-1}(b_p + \hat{F}_p^{\max}) \\ &= M^T A_p^{-1}(b_p + F_p^{\max} + \theta) \\ &= M^T A_p^{-1}\theta_p + M^T A_p^{-1}(b_p + F_p^{\max}) = P\theta_p + w \end{aligned} \tag{16}$$

where $P = M^T A_p^{-1}, w = M^T A_p^{-1}(b_p + F_p^{\max})$.

As mentioned above $\theta = (\theta_1, \theta_2, ..., \theta_K)^T$, write the equation (15) in the following form

$$z^{*}(\theta)=M^{T}A_{p}^{-1}\theta_{p}+M^{T}A_{p}^{-1}(b_{p}+F_{p}^{\max})$$

$$=\left(C^{T}\ \mathrm{SC}^{T}\right)\begin{pmatrix}A_{p_1}\\ H_{p_2}\end{pmatrix}^{-1}\begin{pmatrix}\theta_{S_1}\\ \theta_{S_2}\\ \vdots\\ \theta_{S_{p_1}}\\ 0\end{pmatrix}+\left(C^{T}\ \mathrm{SC}^{T}\right)\begin{pmatrix}A_{p_1}\\ H_{p_2}\end{pmatrix}^{-1}\begin{pmatrix}b_{p_1}+\hat{F}_{p_1}^{\max}\\ H_{p_2}\hat{D}\end{pmatrix} \tag{17}$$

$$If\ S_i\in\{1,2,..,K\cdot T\},\theta_{S_i}=\theta_{S_i},S'=\{S_i\mid S_i\in\{1,2,..,K\cdot T\}\}\ or\ \theta_{S_i}=0$$

Wherein, the subscript of $\theta$ corresponds to the subscript of the active constraint row number of A. In the same critical region, $A_p$ and $b_p+\hat{F}_p^{\max}$ are same.

$$\nabla_{\theta}z^{*}=\left(\left(C^{T}\ UC^{T}\ \mathrm{SC}^{T}\right)\begin{pmatrix}A_{p_1}\\ H_{p_2}\end{pmatrix}^{-1}\right)^{T}=\left(\begin{pmatrix}A_{p_1}\\ H_{p_2}\end{pmatrix}^{-1}\right)^{T}\left(C^{T}\ UC^{T}\ \mathrm{SC}^{T}\right)^{T} \tag{18}$$

For the k-th theta, $if\ \ S_i\in\{(k-1)*T+1,\cdots,k*T\}$.

$$\begin{aligned}\nabla_{\theta_k}\overline{z}^{*}&=\sum_{S_i=(k-1)*T+1}^{k*T}\left(\begin{pmatrix}A_{p_1}\\ H_{p_2}\end{pmatrix}^{-1}\right)_i^{T}\left(C^{T}\ UC^{T}\ \mathrm{SC}^{T}\right)^{T},S_i\in S'\\ \nabla_{\theta_{other}}\overline{z}^{*}&=0,other\notin S'\end{aligned} \tag{19}$$

By evaluating the value of different variable scalar $\nabla_{\theta_k}\overline{z}^{*}$ and providing maximal social welfare with minimal investment, we may choose which lines' capacity should be improved. Assume that inside parameter space, the active and inactive constraints remain constant (i.e., in the given critical region). In that case, we will modify parameter $\boldsymbol{\theta}_k$ with the highest rate/slope of change of the objective function $z$ under the effect of the parameters since increasing the capacity of that line will result in maximum return on investment. This is yet another contribution and innovation in determining the best candidate transmission line for capacity expansion.

Suppose multiple lines' capacity can be increased simultaneously, using the obtained objective function and the investment cost function regarding the parameters. In that case, we will find the optimal solution in the following form.

$$\begin{aligned}&\min\ z=F(\theta_1,\theta_2\ldots\theta_K)\\ &s.t.\begin{cases}IC_1(\theta_1)+IC_2(\theta_2)+\cdots+IC_k(\theta_K)\le BUD\\ 0\le\theta_k\le\overline{\theta}_k\ \ (k=1,2,\cdots,K)\end{cases}\end{aligned} \tag{20}$$

Where, $IC_k(\theta_K),(k=1,2,\cdots,K)$ represents the investment cost function for line k, $0\le\theta_k\le\overline{\theta}_k$ denotes the potential capacity space for line K that can be increased. The highest budgetary limit for the transmission line expansion planning project is represented by $BUD$.

This is an optimization problem regarding multiple parametric as decision variables. We can get the best solutions for all uncertainty parameters, and then we adjust the transmission capacity of each line by $(\theta_1, \theta_2 \ldots \theta_K)^*$. As a result, the lowest cost of electricity generation is achieved within the limited budget, and line congestion is reduced.

Next, we will derive the formula for the critical region. Based on the simplex method, parametric problems are solved by the optimal basic matrix, optimal solution, optimal value, and corresponding critical region derived by the simplex method since the dual solution does not depend on $\theta$ for RHS parametric linear programming problems. Two optimal basic matrixes are considered neighboring matrices when they are both optimal bases for the same $\theta^*$ in parametric space. One can transform into the other one in one dual step (Gal and Nedoma, 1972).

If B is an optimal basic matrix and $m$ is the related index set of basic variables, then (21), which means that all variables must be non-negative, is the primal condition of $\theta$. In a programming problem, this is also the initial condition for decision variables.

At the same time, the parameters themselves have certain constraints shown in (22), these two sets of constraints form the range of values for parameters $\theta$ under the same set of base variables (Gal and Nedoma, 1972).

$$\omega^* = B_f^{-1} b(\theta) \geq 0 \tag{21}$$

$$G\theta \leq g \tag{22}$$

In order to use the simplex approach, we must first transform linear programming to a standard type through using slackness variables to convert all non-equal constraints into equation constraints.

When employing the simplex method, however, the simplex tableau only exhibits the equation constraints, not the variable's possible range of values. Thus, the inequality constraints of Eq. (21) will display the parameters' possible space by ensuring those decision variables are non-negative. However, in the process of solving with the simplex method, the simplex tableau produced only reflects the equation constraints but does not reflect the range of values of the variables, so inequality Eq. (20) will reflect the range of values for the variable, it indicates that variables relate to basis $B_p$ within such the range of values corresponding to parameters $\theta$, so the constraints (20) and (21) reflect the critical region.

In the Lagrange method, we do not normalize the original inequality constraints, and the inequality constraints contain the feasible space of the decision variables. However, we find that the inactive constraint is indeterminate, and under the same active constraints and different parameters, the value on the left of the inactive constraint will change or even not satisfy the constraint. This reflects that the parameters are within a certain range, and the active constraints are unchanged.

Previous studies addressing multi-parametric linear programming problems did not describe how to construct critical regions in detail. Because the solution is defined as a piecewise affine formula about parameter $\theta$ in multi-parametric linear programming problems. So, we take $\omega^*$ having solved by Eq. (14) to inactive constraints:

$$A_s\omega^* \le b_s + \hat{F}_s^{\max} \tag{23}$$

The parameter range obtained by inequality of equation (21) is the feasible range of this critical region

$$\begin{aligned} CR(\theta) &= \left\{\theta \in \Theta \mid A_s\omega^* \le b_s + \hat{F}_s^{\max} \right\} \\ &= \left\{\theta \in \Theta \mid A_s A_p^{-1}(b_p + F_p^{\max} + \theta_p) \le b_s + F_s^{\max} + \theta_s \right\} \\ &= \left\{\theta \in \Theta \mid A_s A_p^{-1}\ \theta_p - \theta_s \le b_s + F_s^{\max} - A_s A_p^{-1}(b_p + F_p^{\max}\ )\right\} \\ &= \left\{\theta \in \Theta \mid \theta_s{}' \le b_s + F_s^{\max} - A_s A_p^{-1}(b_p + F_p^{\max}\ )\right\} \end{aligned} \tag{24}$$

$$\bigcup CR_m = K, K \in \Xi \tag{25}$$

where $s$ is the index set of inactive constraints in each non-overlapping critical region, $m$ denotes the index of critical regions, and K refers the set of all feasible parameters. $\Xi$ represents the entire parameter space, since problems have no feasible solution in some parameter regions, so, there has $K \in \Xi$.

When the RHS of inactive constraint $A_s\omega^* \le b_s + \hat{F}_s^{\max}$ changes due to parameter $\theta_s$ variation, both active and inactive constraints are altered concurrently; this allows us to obtain new variable solutions and critical regions until we have covered the full parameter space. This is analogous to the simplex method's modification in the basic variables. The basic matrix will change when the value of a basic variable goes negative.

If a new integer solution is discovered in the critical region, as shown in Figure 1, compare and update the upper bounds in this region before removing it from the unknown parameter space. Finally, if no solution is found, the procedure should be terminated; this problem has no feasible solution because the full relaxation problem has no solution, then the mixed-integer problem will definitely be infeasible.

### *1.2 Determine the upper bound*

Next, we can use the optimal solution solved by the full relaxed problem at this node to determine the upper boundary in parameter space with a feasible solution. This step, as the first in our approach to performing parametric analysis for large-scale MILP, made the model both novel and practicable.

1) To obtain the lower boundary of the objective function, solve a relaxed mp-LP problem. We can get critical regions $CRs$, corresponding optimal decision variables $G$, and objective functions $z(\theta)$ expressed in terms of parameters after solving the multi-parametric problem at each node.

2) To get approximate integer decision variables by applying the following rule. For each decision

variable that the relaxed problem yielded, if $U_m^* \geq \xi$ ($0 \leq \xi \leq 1$), we update obtained $U_m^*$ to one in the m-th critical region; otherwise, if $U_m^* < \xi$, we let $U_m^*$ to zero. Build and solve the integer linear programming (ILP) model about variables $V$ with determining $U$, finding the optimal integer solution $V$. Combine the variables U and V that emerge from this issue and term it $UV_m$.

$$\begin{aligned} &\min \ V_{12} + V_{13} + \cdots + V_{1T} + \cdots + V_{NT} \\ &st. \begin{cases} U_{it} - U_{i(t-1)} \leq V_{it}, \ t = 2,3,\cdots,T \\ \sum_{k=t-T_i^{on}+1}^{t} V_{ik} - U_{it} \leq 0, t = T_i^{on}+1, T_i^{on}+2,\cdots,T \\ \sum_{k=t-T_i^{off}+1}^{t} V_{ik} + U_{i(t-T_i^{off})} \leq 1, t = T_i^{off}+1, T_i^{off}+2,\cdots,T \end{cases} \end{aligned} \tag{26}$$

3) Substitute $UV_m$ into the mp-LP problem at this node in the m-th critical region and solve the integer problem in region $CR_m$, we can get corresponding optimal updated decision variable $\overline{G}_m$, and the sub-region upper objective function $\overline{z}_i(\theta)$ with parameters $\theta$ in the critical regions $CR_m$.

4) At node i, the problem's upper boundary $\overline{z}_i(\theta)$ is produced by combining the sub-region upper boundary objective functions $\overline{z}_{i,m}(\theta)$ in all regions $CR_{i,m}$.

***1.3 Compare the objective function value of the current solution with the current upper bound***

Obviously, we can get the lower boundary of the original problem at the root node. Here we set the upper bound $\overline{z}_i(\theta)$ at the root node obtained in step 1.2 as the initial upper bound $\overline{z}(\theta)$ of the problem.

For a particular node i, comparing the objective value $\hat{z}_i(\theta)$ which was obtained by relaxing the problem and the related upper boundary $\overline{z}(\theta)$. To simplify the problem, we select some random points in the parameter space, $Q$ denotes the number of random points, first calculate the upper and lower bound errors of each point, and then compute the final error rate using a weighted average of all selected random points, as shown in Equation (32).

$$\delta = \left( \sum_{j=1}^{Q} \frac{\overline{z}(\theta) - \hat{z}_i(\theta)}{\hat{z}_i(\theta)} \right) \Big/ Q \tag{27}$$

If there is $\delta \leq \alpha$, the iteration will be terminated; otherwise, move on to step 2. In this case, $\alpha$ is an artificial threshold or error tolerance boundary that will be adjusted to reflect actual market objectives or computation time burden, which is an innovation from this study.

***Step2. Select a branching variable and update the set of nodes N***

Choose a binary variable from U and generate two new nodes fixing variable selected to be 0 and 1. Add the new nodes to N, and update $n = n + 2$ (i.e., Add two nodes to a branch).

Set the unexplored parameter space of these nodes, $\Xi_{n+1}$ and $\Xi_{n+2}$, equal to the remaining critical regions since some parameter space has been removed in any node i if one of the following conditions is satisfied:

(i) The problem is not feasible.

(ii) An integer solution is found.

(iii) The solution of the node is greater than the current best upper bound in the same space.

This simplifies the operation steps and is easier than iterating through all the nodes in all parameter space.

***Step3.Determinne at nodes***

If $N=\{\varnothing\}$, stopping the iteration with solutions; otherwise, take node $i$ from $N$ and remove it from the list: $n=n-1$. Here, n is the number of nodes that have not been calculated.

***Step4.Solve and compare***

Bringing the determined binary variables into the original model, the model (33) can be obtained, $\omega^{'}$ denotes the variable vector after the binary variables determined by the substitution, solve the mp-LP problem with unknown parameter space $\Xi_i$ related to node i like in step 1 using Lagrange method, if the problem is infeasible for all $\theta\in\Xi_i$, go back to step 3. If the problem is feasible, determine and update the upper bound like in step 1, the upper bound obtained at node i is recorded as $\bar{z}_i(\theta)$. For $\bar{z}_i(\theta)$ and current upper bound $\bar{z}(\theta)$, they are segmented imitation functions in the region. Merge the function and select a smaller upper bound, eventually get a new current upper bound $\bar{z}(\theta)$ (Dua and Pistikopoulos, 2000). Calculate the error between the optimal function value $\hat{z}_i(\theta)$ at node i and the current upper bound $\bar{z}_i(\theta)$; if $\delta\le\alpha$ is held, the iteration end; otherwise, continue.

$$\hat{z}(\theta)=\min_{G,U,V} M^T\omega' \quad st.\begin{cases} A\omega'\le b+\hat{F}^{\max} \\ H\omega'=H\hat{D} \end{cases} \tag{28}$$

Compare the solution of node i with the current upper boundary $\bar{z}(\theta)$. Suppose the optimal solution of transmission planning problem at node i as $\hat{z}_i(\theta)$, formed by m optimal value linear functions $\hat{z}_{i,m}(\theta)$, their corresponding critical regions, $CR_{i,m}$. Record the current upper bound as $\bar{z}(\theta)$, defined by a set of linear functions, $\bar{z}_{m'}(\theta)$ , and their corresponding critical regions, $CR_{m'}^{UB}$.

$$CR^{\text{int}}=CR_{i,m}\bigcap CR_{m'}^{UB} \tag{29}$$

If $CR^{\text{int}}\ne\varnothing$, it is then to check whether $\hat{z}_{i,m}(\theta)\le\bar{z}_{m'}(\theta)$. This can be reflected by a new constraint as following, if the solution is greater than the upper bound in the m-th critical region, remove it from uncertain space.

$$\hat{z}_{i,m}(\theta) - \bar{z}_{m'}(\theta) \le 0 \tag{30}$$

***Step5. Determine nodes and critical region of $\theta$***

If all values of binary variables $U$ are fixed at node *i* we will update the upper bound using the comparison method solution from Step 4. Return to step 3 and jump to an uncalculated node if all critical regions of node *i* are eliminated; otherwise, return to step 2 and build a new branch.

## 5. Case study

The proposed methodologies and results are validated in Section 5.1 by using a single period UCED case to portray the calculation procedure in detail through the synthesis data, which has two buses and two uncertainty transmission lines. The IEEE-5 bus UCED case with 24 periods is also used in Section 5.2 to demonstrate the results and insights; for simplicity, we only consider one line's capacity uncertainty.

### 5.1 Single period synthesis case

Consider the following MILP parametric analysis example:

$$z(\theta) = \min_{x,y} 3x_1 + 5x_2 + 18y_1 + 15y_2$$

$$s.t. \quad \begin{cases} x_1 + x_2 = 15 = D_1(7) + D_2(8) \\ 0.8x_1 + 0.7x_2 \le 11.2 + \theta_1 \\ 0.6x_1 + 0.9x_2 \le 12.4 + \theta_2 \\ x_1 - 10y_1 \le 0 \\ x_2 - 10y_2 \le 0 \\ y_i \in \{0,1\}, i = 1,\ 2 \\ 0 \le \theta_s \le 10,\ s = 1,\ 2 \\ x_j \ge 0, j = 1,\ 2 \end{cases} \tag{31}$$

In this single period problem, we shall ignore the constraints for start-up ramp rate, shut-down ramp rate, and the relationship between UC of generators and start-up state in a single period problem, removing the contrasts for minimum on time and minimum off time. So, we only consider the energy balance and generation constraints. In this synthesis case, two generators and two buses are connected through two lines. To show our methods proposed in section 4, we incorporate two uncertainty parameters on the RHS of the transmission lines' capacity constraints to illustrate the uncertainty of lines' capacity. Both parameters ranged between 0 and 10, indicating that the line's capacity can be enhanced between 0 and 10.

We will solve this problem based on the following steps:

**Step0**: Initialize the collection of nodes that need to be solved, $N = \{\varnothing\}$, $n = 0$.

**Step1**: Solve the fully relaxed problem as follows (at node 0).

$$z(\theta) = \min_{x,y} 3x_1 + 5x_2 + 18y_1 + 15y_2$$

$$s.t. \begin{cases} x_1 + x_2 = 15 = D_1(7) + D_2(8) \\ 0.8x_1 + 0.7x_2 \le 11.2 + \theta_1 \\ 0.6x_1 + 0.9x_2 \le 12.4 + \theta_2 \\ x_1 - 10y_1 \le 0 \\ x_2 - 10y_2 \le 0 \\ -y_1 \le 0 \\ -y_2 \le 0 \\ y_1 \le 1 \\ y_2 \le 1 \\ 0 \le \theta_s \le 10,\ s = 1,\ 2;\ x_n \ge 0,\ n = 1,\ 2 \end{cases} \tag{32}$$

Based on what we introduced in Section 4, we rewrote the format of the above problem as follows. First, the vector of objective function coefficients and the decision variable are:

$$M^T = (3 \quad 5 \quad 18 \quad 15); \omega^T = (x_1 \quad x_2 \quad y_1 \quad y_2)$$

The matrix in the left side of the constant is

$$A = \begin{pmatrix} 0.8 & 0.7 & 0 & 0 \\ 0.6 & 0.9 & 0 & 0 \\ -1 & 0 & 0 & 0 \\ 0 & -1 & 0 & 0 \\ 1 & 0 & -10 & 0 \\ 0 & 1 & 0 & -10 \\ 0 & 0 & -1 & 0 \\ 0 & 0 & 0 & -1 \\ 0 & 0 & 1 & 0 \\ 0 & 0 & 0 & 1 \end{pmatrix}$$

Demand on each bus and transmission line limit corresponding to parameters are:

$$D = \begin{pmatrix} 7 \\ 8 \end{pmatrix}; \hat{F}^{\max} = \begin{pmatrix} 0 + \theta_1 \\ 1 + \theta_2 \end{pmatrix}$$

Based on the Eq. (3), multiply the variable by the matrix $A$ to get

$$A\omega = \begin{pmatrix} 0.8 & 0.7 & 0 & 0 \\ 0.6 & 0.9 & 0 & 0 \\ -1 & 0 & 0 & 0 \\ 0 & -1 & 0 & 0 \\ 1 & 0 & -10 & 0 \\ 0 & 1 & 0 & -10 \\ 0 & 0 & -1 & 0 \\ 0 & 0 & 0 & -1 \\ 0 & 0 & 1 & 0 \\ 0 & 0 & 0 & 1 \end{pmatrix} \begin{pmatrix} x_1 \\ x_2 \\ y_1 \\ y_2 \end{pmatrix} = \begin{pmatrix} 0.8x_1 + 0.7x_2 \\ 0.6x_1 + 0.9x_2 \\ -x_1 \\ -x_2 \\ x_1 - 10y_1 \\ x_2 - 10y_2 \\ -y_1 \\ -y_2 \\ y_1 \\ y_2 \end{pmatrix} \leq \begin{pmatrix} 11.2 + \theta_1 \\ 12.4 + \theta_2 \\ 0 \\ 0 \\ 0 \\ 0 \\ 0 \\ 0 \\ 1 \\ 1 \end{pmatrix} = b + \hat{F}^{\max} \tag{33}$$

By moving the right-hand constant of the inequalities and equality to the left based on the format of Eq. (3), we get the following:

$$\begin{cases} f(\omega) = 3x_1 + 5x_2 + 18y_1 + 15y_2 \\ g(\omega) = \begin{pmatrix} 0.8x_1 + 0.7x_2 - 11.2 - \theta_1 \\ 0.6x_1 + 0.9x_2 - 12.4 - \theta_2 \\ -x_1 \\ -x_2 \\ x_1 - 10y_1 \\ x_2 - 10y_2 \\ -y_1 \\ -y_2 \\ y_1 - 1 \\ y_2 - 1 \end{pmatrix} \leq 0 \\ h(\omega) = x_1 + x_2 - 15 = 0 \end{cases} \tag{34}$$

Next, we will get the following Lagrange function based on Eq. (9).

$$\begin{aligned} L = {} & 3x_1 + 5x_2 + 18y_1 + 15y_2 + \mu_1(0.8x_1 + 0.7x_2 - 11.2 - \theta_1) + \mu_2(0.6x_1 + 0.9x_2 - 12.4 - \theta_2) - \mu_3 x_1 - \mu_4 x_2 \\ & + \mu_5(x_1 - 10y_1) + \mu_6(x_2 - 10y_2) + \mu_7(-y_1) + \mu_8(-y_2) + \mu_9(y_1 - 1) + \mu_{10}(y_2 - 1) + \lambda(x_1 + x_2 - 15) \end{aligned} \tag{35}$$

Then, we can obtain the first derivative of the objective function, inequality and equation with respect to the variable $\omega$ in the original model through Eq. (10).

$$\begin{cases} f(\omega^*) = M^T \omega^*, \nabla f(\omega^*) = \begin{pmatrix} 3 & 5 & 18 & 15 \end{pmatrix}^T \\ g(\omega^*) = A\omega^* - b - \hat{F}^{\max}, \nabla g(\omega^*) = \begin{pmatrix} 0.8 & 0.6 & -1 & -1 & 1 & 0 & 0 & 0 & 0 & 0 \\ 0.7 & 0.9 & 0 & 0 & 0 & 1 & 0 & 0 & 0 & 0 \\ 0 & 0 & 0 & 0 & -10 & 0 & -1 & 0 & 1 & 0 \\ 0 & 0 & 0 & 0 & 0 & -10 & 0 & -1 & 0 & 1 \end{pmatrix} \\ h(\omega^*) = H\omega^* - H\hat{D}, \nabla h(\omega^*) = \begin{pmatrix} 1 & 1 & 0 & 0 \end{pmatrix}^T \end{cases} \tag{36}$$

According to the Lagrange equation, the following KKT conditions are derived by Eqs. (11a)-(11e).

$$\begin{cases} \nabla_x L = \dfrac{\partial L}{\partial \omega} = 0 \Rightarrow M + A^T \mu + H^T \lambda = 0 \\ \nabla_{\lambda_j} L = \dfrac{\partial L}{\partial \lambda_j} = h_j(\omega) = H\omega - H\hat{D} = 0, j = 1 \\ g_i(\omega) = A\omega - b - \hat{F}^{\max} \le 0, i = 1, ..., K + 4N \\ \mu_i \ge 0, i = 1, ..., K + 4N \\ \mu_i g_i(\omega) = 0 \Leftrightarrow \mu^T \left( A\omega - b - \hat{F}^{\max} \right) = 0, i = 1, ..., K + 4N \end{cases} \tag{37}$$

When $\theta_1 = 0, \theta_2 = 0$, we can obtain the solutions of the variables and the Lagrange multipliers corresponding to the equation constraint and the inequality constraint. Also, we can get the active and inactive constraints by the value of the Lagrange multiplier.

$$\begin{cases} x_1 = 7 & \mu_1 = 9 \\ x_2 = 8 & \mu_2 = 0 \\ y_1 = 0.7 & \mu_5 = 2.4 \\ y_2 = 0.8 & \mu_6 = 1.5 \\ \lambda = -12.8 & \mu_3 = \mu_4 = \mu_7 = \mu_8 = \mu_9 = \mu_{10} = 0 \end{cases} \tag{38}$$

Meanwhile, we can extract the active constraints from Eq. (38) using the Lagrange multipliers: the first, fifth, and sixth inequalities, as well as the unique equality constraint coefficients, form the basis matrix in Eqs. (13) and (14) because the Lagrange multipliers for these constraints are not zero.

$$A_p = \begin{pmatrix} 0.8 & 0.7 & 0 & 0 \\ 1 & 0 & -10 & 0 \\ 0 & 1 & 0 & -10 \\ 1 & 1 & 0 & 0 \end{pmatrix}, b_p + F_p^{\max} + \theta_p = \begin{pmatrix} 11.2 + \theta_1 \\ 0 \\ 0 \\ 15 \end{pmatrix}, \theta_p = \begin{Bmatrix} \theta_{S_1} \\ \theta_{S_2} \\ \theta_{S_3} \\ 0 \end{Bmatrix} = \begin{Bmatrix} \theta_1 \\ \theta_5 \\ \theta_6 \\ 0 \end{Bmatrix} = \begin{Bmatrix} \theta_1 \\ 0 \\ 0 \\ 0 \end{Bmatrix} \tag{39}$$

Similarly, we can also get inactive constraints, which correspond to the matrix composed of the coefficients of the second to fourth and seventh to tenth inequality constraints from Eq. (38). We also computed the inverse matrix $A_p$ for subsequent operations at the same time.

$$A_s = \begin{pmatrix} 0.6 & 0.9 & 0 & 0 \\ -1 & 0 & 0 & 0 \\ 0 & -1 & 0 & 0 \\ 0 & 0 & -1 & 0 \\ 0 & 0 & 0 & -1 \\ 0 & 0 & 1 & 0 \\ 0 & 0 & 0 & 1 \end{pmatrix}, b_s + F_s^{\max} + \theta_s = \begin{pmatrix} 12.4 + \theta_2 \\ 0 \\ 0 \\ 0 \\ 0 \\ 1 \\ 1 \end{pmatrix}, A_p^{-1} = \begin{pmatrix} 10 & 0 & 0 & -7 \\ -10 & 0 & 0 & 8 \\ 1 & -0.1 & 0 & -0.7 \\ -1 & 0 & -0.1 & 0.8 \end{pmatrix} \tag{40}$$

Using Eq. (15) in Section 4, we will derive the best decision variable $\omega_{0,1}^*$, as shown in Eq. (39).

$$\omega_{0,1}^{*} = A_p^{-1}(b_p + F_p^{\max} + \theta_p) = \begin{pmatrix} 10\cdot\theta_1 + 7 \\ 8 - 10\cdot\theta_1 \\ \theta_1 + 7/10 \\ 4/5 - \theta_1 \end{pmatrix} \tag{41}$$

Based on Eq. (24), the critical region is obtained as follows:

$$CR_{0,1} = \left\{\theta \in \Theta \mid A_s\omega^* \le b_s + \hat{F}_s^{\max}\right\} = \left\{\theta \in \Theta \,\middle|\, \begin{pmatrix} 57/5 - 3\cdot\theta_1 \\ -10\cdot\theta_1 - 7 \\ 10\cdot\theta_1 - 8 \\ -\theta_1 - 7/10 \\ \theta_1 - 4/5 \\ \theta_1 + 7/10 \\ 4/5 - \theta_1 \end{pmatrix} \le \begin{pmatrix} 12.4 + \theta_2 \\ 0 \\ 0 \\ 0 \\ 0 \\ 1 \\ 1 \end{pmatrix}\right\} = \left\{0 \le \theta_1 \le 0.3, 0 \le \theta_2 \le 10\right\} \tag{42}$$

Based on the optimal solutions in Eq. (39) and the critical region in Eq. (40), the optimal objective function is shown as

$$z_{0,1} = 85.6 - 17.0\cdot\theta_1, \mathit{when}\left\{0 \le \theta_1 \le 0.3, 0 \le \theta_2 \le 10\right\} \tag{43}$$

Here, we recognize the theta as a decision variable and then calculate the first-order response function of $z_{0,1} = 85.6 - 17.0\cdot\theta_1$ on theta through Eq. (19). Then, we find the cost reduction rate with the upgrade of the transmission line's capacity in the same critical region.

Here, $S$ represents the set of active constraints from $A$ in the constraints. From Eq (39), we will get the following results

$$S = \left\{S_1, S_2, S_3\right\} = \left\{1,5,6\right\} \tag{44}$$

Here, $S'$ denotes inequality constraint row subscripts with parameters on the right-hand side.

$$S' = \left\{S_i \mid S_i \in \left\{1,2\right\}\right\} \tag{45}$$

For the first parameter, $\theta_1$, $S_1 = 1 \in S'$.

$$\nabla_{\theta_1}\bar{z}^* = \left(\begin{pmatrix} A_{p-1} \\ H \end{pmatrix}^{-1}\right)_1^T \left(C^T\ \mathrm{SC}^T\right)^T = -17, S_1 \in S' \tag{46}$$

For the second parameter, $\theta_2$, because $\forall\, S_i \notin S'$.

$$\nabla_{\theta_2}\bar{z}^* = 0 \tag{47}$$

The results (i.e., Eq. (46) and Eq. (47)) reveal that increasing the capacity of the first line by one unit reduces electricity generation costs by 17 units, however, increasing the capacity of the second line has no effect on system costs.

When parameter $\theta_1 = 0.3$, the fourth inactive constraint in matrix $A_s$ of Eq. (39), becomes an active constraint, the first active constraint in matrix $A_p$ of Eq. (39) becomes an inactive constraint, and the matrix $A_p$ in Eq. (39) changes with it. Following, the updated optimal Lagrange multipliers are shown as follows:

$$\begin{cases} x_1 = 10 & \mu_2 = 0 \\ x_2 = 5 & \mu_5 = 3.5 \\ y_1 = 1 & \mu_6 = 1.5 \\ y_2 = 0.5 & \mu_7 = 11 \\ \lambda = -6.5 & \mu_1 = \mu_3 = \mu_4 = \mu_7 = \mu_8 = \mu_9 = \mu_{10} = 0 \end{cases} \tag{48}$$

At this time (i.e., $\theta_1 = 0.3$), the revised matrixes obtained based on the active constraints and inactive constraints are shown as follows, respectively. According to the Lagrange multipliers obtained from Eq. (48), we can obtain the active constraints as follows, the fifth, sixth and seventh inequalities and unique equality constraint coefficients constitute the basis matrix by Eq. (13) and (14).

$$A_p = \begin{pmatrix} 1 & 0 & -10 & 0 \\ 0 & 1 & 0 & -10 \\ 0 & 0 & 1 & 0 \\ 1 & 1 & 0 & 0 \end{pmatrix}, b_p + F_p^{\max} + \theta_p = \begin{pmatrix} 0 \\ 0 \\ 1 \\ 15 \end{pmatrix}, \theta_p = \begin{Bmatrix} \theta_{S_1} \\ \theta_{S_2} \\ \theta_{S_3} \\ 0 \end{Bmatrix} = \begin{Bmatrix} \theta_5 \\ \theta_6 \\ \theta_7 \\ 0 \end{Bmatrix} \tag{49}$$

At the same time, we can also get inactive constraints, corresponding to the matrix consisting of the coefficients of the first to fourth and eighth to tenth inequality constraints from Eq. (48). We also calculated the inverse matrix $A_p$.

$$A_s = \begin{pmatrix} 0.8 & 0.7 & 0 & 0 \\ 0.6 & 0.9 & 0 & 0 \\ -1 & 0 & 0 & 0 \\ 0 & -1 & 0 & 0 \\ 0 & 0 & -1 & 0 \\ 0 & 0 & 0 & -1 \\ 0 & 0 & 0 & 1 \end{pmatrix}, b_s + F_s^{\max} + \theta_s = \begin{pmatrix} 11.2 + \theta_1 \\ 12.4 + \theta_2 \\ 0 \\ 0 \\ 0 \\ 0 \\ 1 \end{pmatrix}, A_p^{-1} = \begin{pmatrix} 1 & 0 & 10 & 0 \\ -1 & 0 & -10 & 1 \\ 0 & 0 & 1 & 0 \\ -0.1 & -0.1 & -1 & 0.1 \end{pmatrix} \tag{50}$$

Similarly, using Eqs. (15), (16), and (24), we will directly calculate the best expression for decision variables $\omega_{0,2}^*$, the critical region $CR_{0,2}$, and the optimal objective value $z_{0,2}$, respectively.

$$\omega_{0,2}^* = A_p^{-1}(b_p + F_p^{\max} + \theta_p) = \begin{pmatrix} 10 \\ 5 \\ 1 \\ 0.5 \end{pmatrix} \tag{51}$$

$$CR_{0,2}=\left\{\theta\in\Theta\,|\,A_s\omega^*\le b_s+\hat{F}_s^{\max}\right\}=\left\{\theta\in\Theta\,\middle|\begin{pmatrix}11.5\\10.5\\-10\\-5\\-1\\-0.5\\0.5\end{pmatrix}\le\begin{pmatrix}11.2+\theta_1\\12.4+\theta_2\\0\\0\\0\\0\\1\end{pmatrix}\right\}=\left\{0.3\le\theta_1\le10,0\le\theta_2\le10\right\} \tag{52}$$

$$z_{0,2}=80.5, when\ \left\{0.3\le\theta_1\le10,0\le\theta_2\le10\right\} \tag{53}$$

In this critical region (i.e., $0.3\le\theta_1\le10,0\le\theta_2\le10$ ), we will investigate how theta/transmission line capacity changes affects the system cost (i.e., optimal objective value).

Because the parameter has not changed, $S$ is known from Eq. (49), and $S'$ is the same as Eq. (45).

$$S=\left\{S_1,S_2,S_3\right\}=\left\{5,6,7\right\} \tag{54}$$

For the first parameter, $\theta_1$, and the second parameter, $\theta_2$, because there is $\forall\ S_i\not\subset S'$, so we will get

$$\nabla_{\theta_1}\overline{z}^*=0 \tag{55}$$

$$\nabla_{\theta_2}\overline{z}^*=0 \tag{56}$$

According to the results in Eqs. (56) and (57), enhancing the capacity of the first and second transmission lines has no effect on optimal electricity generation cost outcome $\overline{z}^*$ in the second critical region (i.e., $0.3\le\theta_1\le10,0\le\theta_2\le10$).

**Step2**: Choosing a node and fixing UC, followed by solving multi-parametric LP as in step 1.

***Node1*** (Fix $y_1=0$)

The problem is infeasible, so this node is ignored.

***Node2*** (Fix $y_1=1$)

$$z(\theta)=\min_{x,y}3x_1+5x_2+18+15y_2$$

$$s.t.\ \begin{cases}x_1+x_2=15=D_1(7)+D_2(8)\\0.8x_1+0.7x_2\le11.2+\theta_1\\0.6x_1+0.9x_2\le12.4+\theta_2\\x_1-10\le0\\x_2-10y_2\le0\\-y_2\le0\\y_2\le1\\0\le\theta_s\le10,\ s=1,\ 2;x_n\ge0,\ n=1,\ 2\end{cases} \tag{57}$$

According to the Lagrange method, we can get the optimal decision variable from Eq. (15), the optimal value derived from Eq. (16) and the critical region based on Eq. (24). The optimal solutions are displayed.

$$\begin{cases} x_1 = 10\cdot\theta_1 + 7, x_2 = 8 - 10\cdot\theta_1, y_1 = 1, y_2 = 0.8 - \theta_1 \\ z_{2,1}(\theta) = 91.0 - 35\cdot\theta_1 \qquad CR_{2,1} = \{\Theta \mid 0 \le \theta_1 \le 0.3, 0 \le \theta_2 \le 10\} \\ x_1 = 10, x_2 = 5, y_1 = 1, y_2 = 0.5 \\ z_{2,2}(\theta) = 80.5 \qquad CR_{2,2} = \{\Theta \mid 0.3 \le \theta_1 \le 10, 0 \le \theta_2 \le 10\} \end{cases} \tag{58}$$

***Node3*** (Fix $y_1 = 1, y_2 = 0$)

The problem is infeasible.

***Node4*** (Fix $y_1 = 1, y_2 = 1$)

$$z(\theta) = \min_{x,y} 3x_1 + 5x_2 + 18 + 15$$

$$s.t. \quad \begin{cases} x_1 + x_2 = 15 = D_1(7) + D_2(8) \\ 0.8x_1 + 0.7x_2 \le 11.2 + \theta_1 \\ 0.6x_1 + 0.9x_2 \le 12.4 + \theta_2 \\ x_1 - 10 \le 0 \\ x_2 - 10 \le 0 \\ 0 \le \theta_s \le 10, \ s = 1, 2 \\ x_n \ge 0, \ n = 1, 2 \end{cases} \tag{59}$$

Similarly, by using Eq. (15) and Eq. (16), the optimal solutions are shown as follows:

$$\begin{cases} x_1 = 10\cdot\theta_1 + 7, x_2 = 8 - 10\cdot\theta_1, y_1 = 1, y_2 = 1 \\ z_{4,1}(\theta) = 94.0 - 20.0\cdot\theta_1 \qquad CR_{4,1} = \{\Theta \mid 0 \le \theta_1 \le 0.3, 0 \le \theta_2 \le 10\} \\ x_1 = 10, x_2 = 5, y_1 = 1, y_2 = 1 \\ z_{4,2}(\theta) = 88.0 \qquad CR_{4,2} = \{\Theta \mid 0.3 \le \theta_1 \le 10, 0 \le \theta_2 \le 10\} \end{cases} \tag{60}$$

Based on the results of Eq. (61), we know that the decision-maker/ ISO should increase the 0.3-unit capacity of the first line. The results of Eq. (61) show that 1) increasing the capacity of the second line has no effect on the reduction of electricity generation costs, indicating that this line is not congested; and 2) in the current situation, increasing the capacity of the first line by more than 0.3 capacity has no effect on the optimization results.

**Step3**: To sum up, we will draw the following optimal three-dimensional plane curve through MPT (see https://www.mpt3.org/Main/HomePage for details).

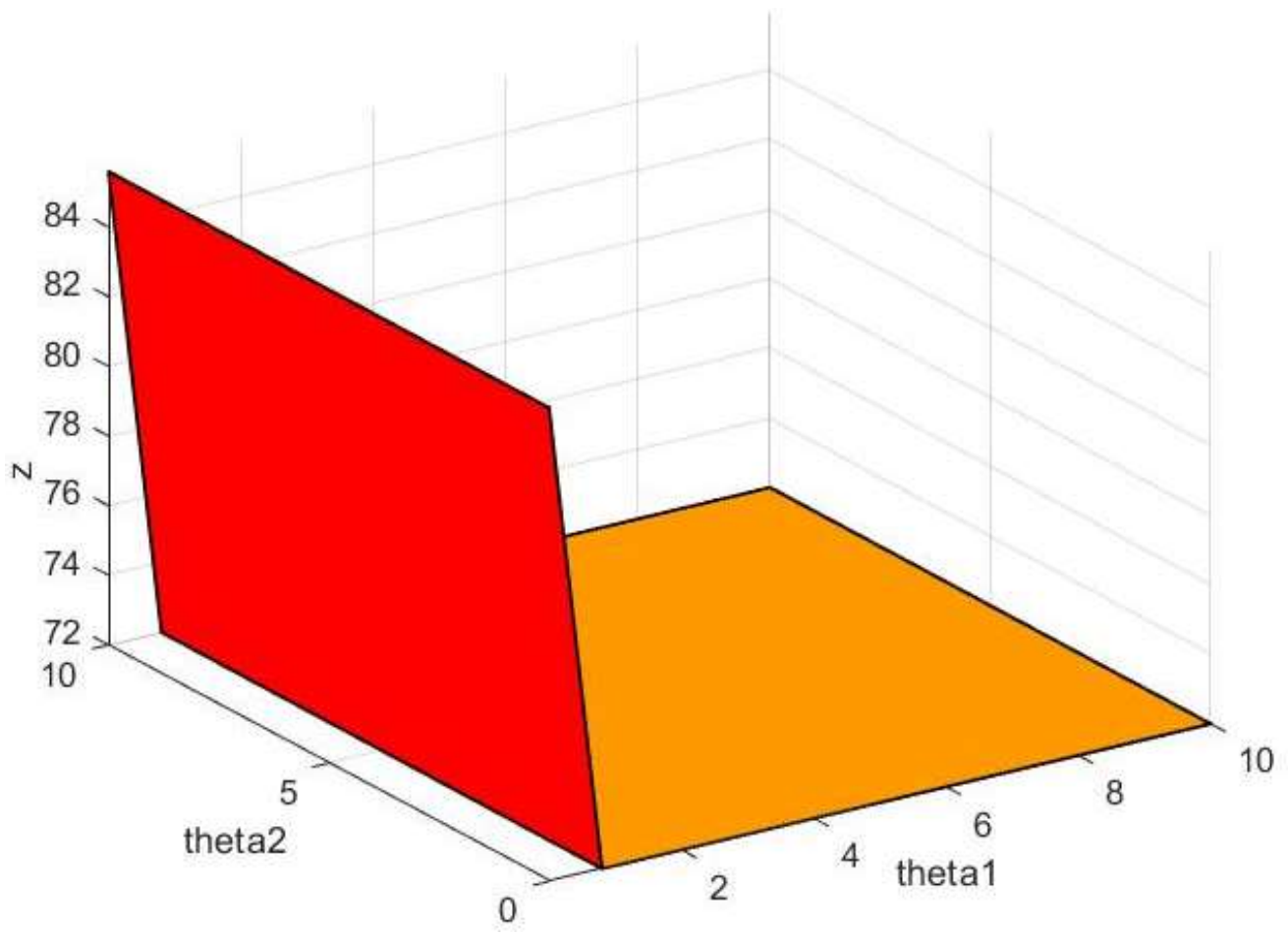

**Figure 2: Optimal solution with parametric analysis**
(at node 0)

$$\begin{cases} z_1 = 85.6 - 17.0 \times \theta_1, where,\ 0 \le \theta_1 \le 0.3, 0 \le \theta_2 \le 10 \\ z_2 = 80.5, where,\ 0.3 \le \theta_1 \le 10, 0 \le \theta_2 \le 10 \end{cases}$$

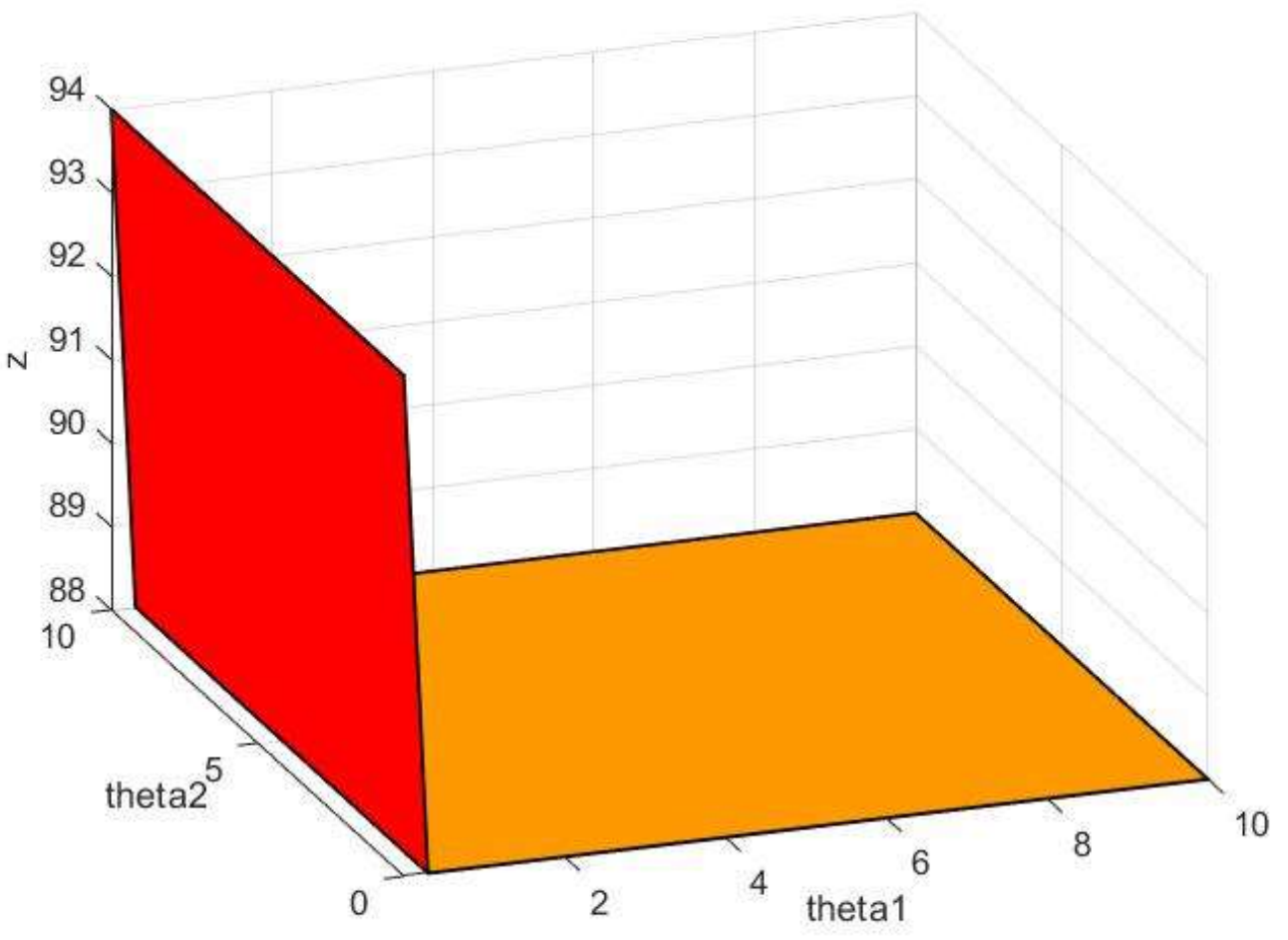

**Figure 3: Optimal solution with parametric analysis**
(at node 4)

$$\begin{cases} z_1 = 94.0 - 20.0 \times \theta_1, where,\ 0 \le \theta_1 \le 0.3, 0 \le \theta_2 \le 10 \\ z_2 = 88.0, where,\ 0.3 \le \theta_1 \le 10, 0 \le \theta_2 \le 10 \end{cases}$$

Figure 2 is the optimal parametric analysis solution considering the uncertainty of two transmission lines' capacity at the root node by relaxing all binary variables to continuous variables. This is a three-dimensional plane, with the three dimensions representing theta1, theta2, and z, respectively. The accurate optimal solution can be observed in this figure. Figure 3 depicts the relation between the parameters ($\theta_1, \theta_2$) and the ideal objective function value at node 4 (i.e., $y_1 = 1, y_2 = 1$). This is the final solution to the mp-MILP problem in Eq. (31), and the obtained solution is also the optimal solution satisfying all constraints. The z-axis coordinate vertex representing the optimal value is higher than the coordinate vertex in Figure 2, indicating that the cost at node 4 is higher than the cost obtained by the root node, which is reasonable because the binary variables are relaxed at the root node.

### 5.2 IEEE-5 bus 24-periods UCED case

This section applies our method to a 24-periods and IEEE-5 bus UCED problem. This problem has five buses and six lines, and all data (i.e., transmission capacity, generation shift factor, per-unit electricity generation cost, etc.) is generated in MATLAB using MATPOWER. To simplify and more intuitively show our algorithm and conclusion, we only consider the uncertainty of the first line transmission capacity (i.e., parametric analysis MILP problem for one parameter). By relaxing all binary variables to continuous variables at the root node between 0 and 1. As a result, we can derive the lower boundary in the entire parameter space based on linear programming parametric analysis.

$$\underline{z}=\begin{cases}461102.11-68.684373\cdot\theta_4, 0\le\theta_4\le 26.145\\ 460653.52-51.526754\cdot\theta_4, 26.145\le\theta_4\le 26.865\\ 460192.22-34.355661\cdot\theta_4, 26.865\le\theta_4\le 31.9985\\ 459642.34-17.171093\cdot\theta_4, 31.9985\le\theta_4\le 44.2195\\ 458883.04, 44.2195\le\theta_4\le 100\end{cases}\tag{61}$$

The parameter space is categorized into five critical regions, and the optimal parametric analysis objective function is a piecewise affine function and split into five sub-piecewise parts. Eq. (61) reveals that the optimal objective function is a constant value when $44.2195\le\theta_4\le 100$.

After the fully relaxed solution is obtained at the root node, employing Eq. (26) and transforming the obtained optimal continuous variables to 0 or 1. Here, we assume $\xi=0$, then fixing these variables, we get the upper bound of the IEEE-5 bus UCED problem.

$$\overline{z}=\begin{cases}461131.17-68.720651\cdot\theta_4, 0\le\theta_4\le 26.145\\ 460681.98-51.540488\cdot\theta_4, 26.145\le\theta_4\le 26.865\\ 460218.44-34.360325\cdot\theta_4, 26.865\le\theta_4\le 31.9985\\ 459668.7-17.180163\cdot\theta_4, 31.9985\le\theta_4\le 44.2195\\ 458909, 44.2195\le\theta_4\le 100\end{cases}\tag{62}$$

When the upper and lower bound at the root node are obtained, we can calculate the error value in each critical region between $\underline{z}$ and $\overline{z}$. Because the current capacity of line 4 is 150MW, the first line segment in Eq. (62) shows that if the capacity of line 4 may be regulated between 150 and 176.145 (i.e., 150+26.145), raising 1(MW) capacity of line 4 will decrease 68.720651 ($). The third line segment indicates that by increasing line 4 capacity from 176.865 (i.e., 150+26.865) to 181.9985 (i.e., 150+31.9985), the system cost will be reduced by 34.360235($). The fifth line demonstrates that increasing the capacity of line 4 would not benefit the system if the existing capacity is already more than 194.2195(i.e.,150+44.2195=194.2195).

Obviously, with the increase of line's capacity marginal effect will decrease through increasing line capacity will reduce the system's cost. As a result, transmission planners or ISO must balance the marginal benefits of increased line capacity against the cost of improving line capacity. Decision makers should increase capacity only if the rate of system cost decline is greater than the cost of increasing per unit capacity.

According to the actual situation, we set $\alpha=0.01$ as the upper limit of the tolerance error, and the error can be calculated using Eq. (27) as follows. In this instance, we compute the error by randomly selecting 50 values (i.e., $Q=50$ ) from each critical region.

$$\delta = \begin{cases} 0.0063\%, 0 \le \theta_1 \le 26.145 \\ 0.00612\%, 26.145 \le \theta_1 \le 26.865 \\ 0.00568\%, 26.865 \le \theta_1 \le 31.9985 \\ 0.00566\%, 31.9985 \le \theta_1 \le 44.2195 \\ 0.00566\%, 44.2195 \le \theta_1 \le 100 \end{cases} \tag{53}$$

Because the error in each crucial zone is less than 0.006 %, which is within the tolerated range (i.e., $\delta < \alpha = 0.01$), the iteration operation will be terminated as specified in step 1 of Section 4. So, compared to traditional B&B methods for exploring all potential nodes, we can conclude that the approximate optimal solution to this problem was obtained quickly based on the proposed method. Furthermore, the accuracy degree of approximation was very high, about 99.994 % for this problem.

Therefore, by using the MPT, we will draw the final optimal approximate results, which is the upper bound curve as shown in Figure 4.

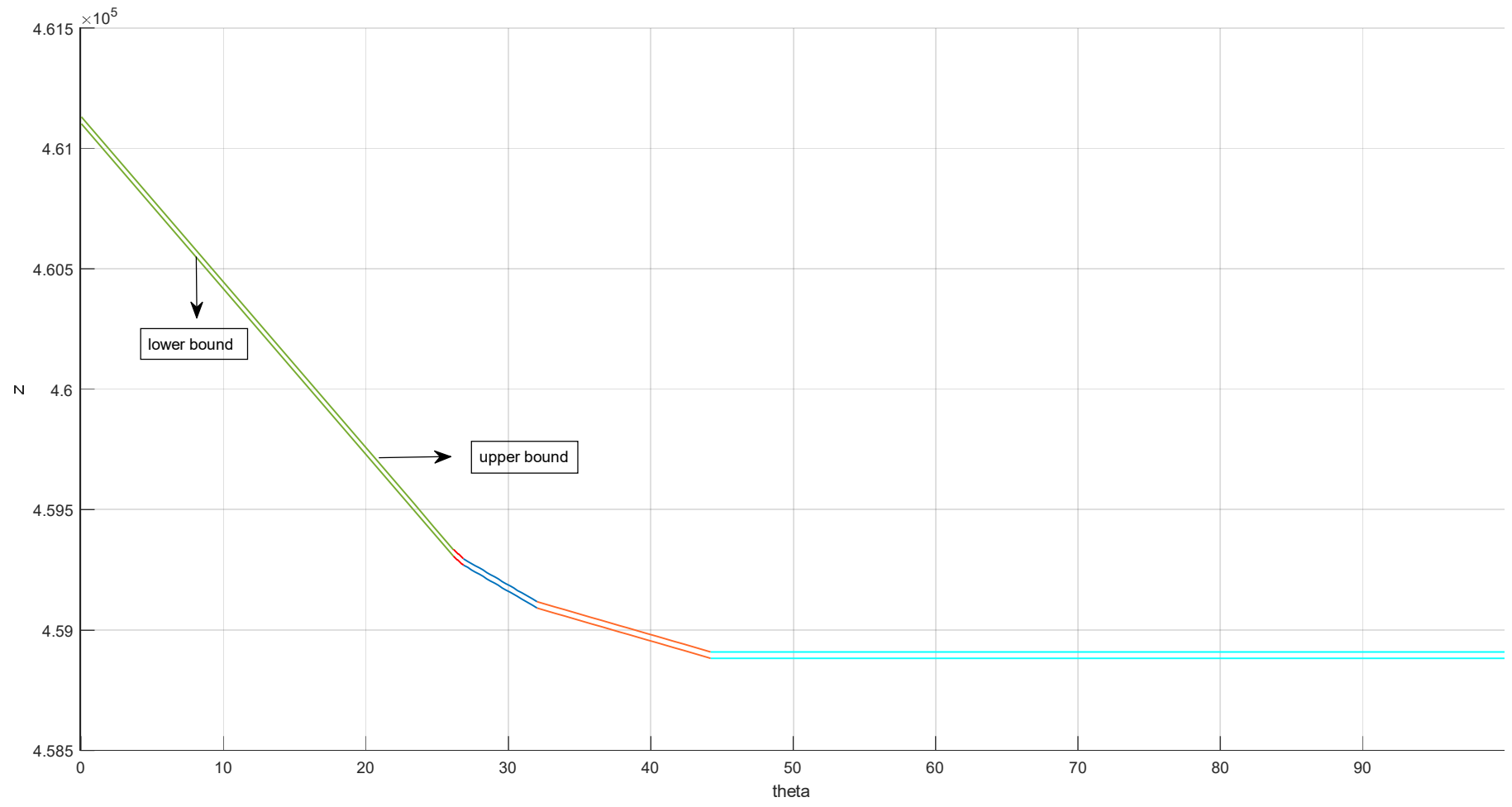

**Figure 4: The lower and upper boundaries at root node**

The horizontal axis of the graph indicates the space of parameter theta, or the possible space for improving transmission line capacity, and the vertical axis represents the best objective value. The difference between the lower and upper bound curves represents the error from Eq. (53), and it can be seen that the distance between the upper and lower bounds is quite modest. The UCED example in this IEEE-5 bus 24-period includes 120 binary decision variables; we should evaluate enormous nodes (i.e., 2.6585e+36) using the classical B&B approach, which brings the simulation burden and cannot be supported in MATLAB using the MPT solver; the optimal solution will fall into the lower and upper bound curves in Figure 4. Based on our approach, only investigating one node yields the appropriate optimal solution, the upper bound curve.

## 6. Conclusion and Future work

This is the first study to use parametric analysis for multiple periods MILP to model transmission expansion planning challenges based on ISO multiple periods UCED problems. To make it more realistic, considering the transmission line capacity uncertainty problem, which also includes starting ramp rate restrictions, ramp-up rate limits, ramp-down limits, and shutdown constraints. This paper presented a novel approach for performing large-scale MILP parametric analysis using the B&B method to streamline the solution process by extracting an approximate upper boundary near the optimum output. The optimal decision structure theorized the classic conclusions of Acevedo and Pistikopoulos (1997), Li and Bo (2009) and Oberdieck et.al (2014), and it differed from conventional policies, referred to as optimal in the literature, in that it did not take into account the enormous binary variables.

Increasing existing transmission line capacity is a desirable technique for responding to development by countering instability and ensuring transmission line security. The optimum upgrade policy for a transmission line planning problem that includes raising the capacity of existing lines corresponds to ISO's perspective on which lines' capacity should be increased to reduce system cost. We characterize this as a RHS uncertainty parametric analysis for MILP. We provide a new matrix approach that employs the classical Lagrangian function and the branch and bound algorithm. We proposed a strategy for achieving the relevant optimal solution and critical regions by employing optimal Lagrange multiples to specify the original matrix's active and inactive constraints.

We extended the classic B&B approach to produce the approximate analytical solution for the large-scale MILP parametric analysis by comparing the error between the lower and upper boundaries at each node. By modifying the tolerated error range, the transmission planner or ISO obtains the analytically best solution for the small-scale UCED problem and achieves the approximate results for the large-scale UCED problem. In that case, the decision maker will get acceptable best results during the limited time to generate the best transmission capacity expansion decisions.

By analyzing the analytically optimal system cost function and the potential improved transmission capacity and comparing the cost-cutting rate to the capacity increase, decision-makers can adopt realistic line capacity changes and achieve minimum system cost. As a result, transmission planners or the ISO must balance the marginal advantages of additional line capacity against the cost of improving line capacity.

Building new lines or increasing the transmission capacity of existing lines can help to relieve network congestion while maintaining system security and stability. In our model, only the upgrading capacity of the current transmission line is considered. We approach this problem using multiple parametric MILP while only accounting for the uncertainty of the RHS of the constraints. Future studies can address establishing additional lines (TEP problem) by adding new constraints and adjusting the generation shift factor (GSG) in the left-hand constraint matrix. The results of this paper are founded based on the assumption that generators have linear generating costs. However, in many studies, such as Hua and Baldick (2017) and Yu et al. (2020), it is likely that the main structural results and the proposed approach of this paper will still hold after relaxing the sampling assumptions. Therefore, another avenue worth exploring would be modeling the electricity generation cost as a nonlinear function.

**Acknowledgments**

This research is based upon work supported by the US Department of Energy's Office of Energy Efficiency and Renewable Energy (EERE) under the Water Power Technologies Office (Grant no. DE-EE0008781). The views expressed herein do not necessarily represent the views of the U.S. Department of Energy or the United States Government.

We thank Haotian Chen, Hossein MehdipourPicha, and Yang Chen (Missouri University of Science and Technology) for their helpful comments and discussion on previous drafts.

**Appendix A: Proposed Matrix Construction Approach**

$$GSF=\begin{pmatrix} \left.\begin{matrix} \overbrace{GSF_{1-1}\ \ 0\cdots0}^{|T|(i=1)} & \overbrace{GSF_{1-2}\ \ 0\cdots0}^{|T|(i=2)} & \cdots\cdots & \overbrace{GSF_{1-N}\ \ 0\cdots0}^{|T|(i=N)} \\ 0\ \ GSF_{1-1}\cdots0 & 0\ \ GSF_{1-2}\cdots0 & \cdots\cdots & 0\ \ GSF_{1-N}\cdots0 \\ \vdots & \vdots & & \vdots \\ 0\cdots0\ \ GSF_{1-1} & 0\cdots0\ \ GSF_{1-2} & \cdots\cdots & 0\cdots0\ \ GSF_{1-N} \end{matrix}\right\}|T|(k=1) \\ \vdots \qquad\qquad \ddots \qquad\qquad \vdots \\ \left.\begin{matrix} GSF_{K-1}\ \ 0\cdots0 & GSF_{K-2}\ \ 0\cdots0 & \cdots\cdots & GSF_{K-N}\ \ 0\cdots0 \\ 0\ \ GSF_{K-1}\cdots0 & 0\ \ GSF_{K-2}\cdots0 & \cdots\cdots & 0\ \ GSF_{K-N}\cdots0 \\ \vdots & \vdots & & \vdots \\ \underbrace{0\cdots0\ \ GSF_{K-1}}_{|T|(i=1)} & \underbrace{0\cdots0\ \ GSF_{K-2}}_{|T|(i=2)} & \cdots\cdots & \underbrace{0\cdots0\ \ GSF_{K-N}}_{|T|(i=N)} \end{matrix}\right\}|T|(k=K) \end{pmatrix}_{KT\times NT} \tag{A. 1}$$

Here, $GSF$ denotes generation shift factor, and the generation shift factors without change for any given generator i-th on the same line in all periods.

In the following matrixes, $UR$ is the ramp-up limit, $SR$ is the startup ramp limit, $G^{\max}$ is the maximum power generation capacity for generators, $DR$ is the ramp-down limit, and $SD$ is the shutdown limit, respectively.

$$UR=\begin{pmatrix} \left.\begin{matrix} \overbrace{UR_1\ \ 0\ \cdots\ 0}^{|T|}\cdots\cdots\ 0 \\ 0\ \ UR_1\ \cdots\ 0\ \cdots\cdots 0 \\ \vdots \qquad \ddots \qquad \vdots \\ 0\cdots\cdots\cdots UR_1\cdots\cdots\ 0 \end{matrix}\right\}|T|(i=1) \\ \vdots \qquad \ddots \qquad \vdots \\ \left.\begin{matrix} 0\ \cdots\ \ UR_N\cdots\cdots\ 0 \\ 0\ \cdots\ \ 0\ UR_N\cdots 0 \\ \vdots \qquad \ddots \qquad \vdots \\ 0\cdots\cdots\ \ \underbrace{0\cdots\cdots\ UR_N}_{|T|} \end{matrix}\right\}|T|(i=N) \end{pmatrix}_{NT\times NT} \quad SR=\begin{pmatrix} \left.\begin{matrix} \overbrace{SR_1\ \ 0\ \cdots\ 0}^{|T|}\cdots\cdots\ 0 \\ 0\ \ SR_1\ \cdots\ 0\ \cdots\cdots 0 \\ \vdots \qquad \ddots \qquad \vdots \\ 0\cdots\cdots\cdots SR_1\cdots\cdots\ 0 \end{matrix}\right\}|T|(i=1) \\ \vdots \qquad \ddots \qquad \vdots \\ \left.\begin{matrix} 0\ \cdots\ \ SR_N\cdots\cdots\ 0 \\ 0\ \cdots\ \ 0\ SR_N\cdots 0 \\ \vdots \qquad \ddots \qquad \vdots \\ 0\cdots\cdots\ \ \underbrace{0\cdots\cdots\ SR_N}_{|T|} \end{matrix}\right\}|T|(i=N) \end{pmatrix}_{NT\times NT}$$

$$G^{\max}=\begin{pmatrix} \left.\begin{matrix} \overbrace{G_1^{\max}\ \ 0\ \cdots\ 0}^{|T|}\cdots\cdots\ 0 \\ 0\ \ G_1^{\max}\ \cdots\ 0\ \cdots\cdots 0 \\ \vdots \qquad \ddots \qquad \vdots \\ 0\ \cdots\cdots\cdots G_1^{\max}\cdots\cdots\ 0 \end{matrix}\right\}|T|(i=1) \\ \vdots \qquad \ddots \qquad \vdots \\ \left.\begin{matrix} 0\ \cdots\ \ G_N^{\max}\cdots\cdots\ 0 \\ 0\ \cdots\ \ 0\ G_N^{\max}\cdots 0 \\ \vdots \qquad \ddots \qquad \vdots \\ 0\cdots\cdots\ \ \underbrace{0\cdots\cdots\ G_N^{\max}}_{|T|} \end{matrix}\right\}|T|(i=N) \end{pmatrix}_{NT\times NT} \tag{A.2}$$

$$DR=\begin{pmatrix} \left.\begin{matrix} \overbrace{DR_1\ \ 0\ \cdots\ 0}^{|T|}\cdots\cdots\ 0 \\ 0\ \ DR_1\ \cdots\ 0\ \cdots\cdots 0 \\ \vdots \qquad \ddots \qquad \vdots \\ 0\cdots\cdots\cdots DR_1\cdots\cdots\ 0 \end{matrix}\right\}|T|(i=1) \\ \vdots \qquad \ddots \qquad \vdots \\ \left.\begin{matrix} 0\ \cdots\ \ DR_N\cdots\cdots\ 0 \\ 0\ \cdots\ \ 0\ DR_N\cdots 0 \\ \vdots \qquad \ddots \qquad \vdots \\ 0\cdots\cdots\ \ \underbrace{0\cdots\cdots\ DR_N}_{|T|} \end{matrix}\right\}|T|(i=N) \end{pmatrix}_{NT\times NT} \quad SD=\begin{pmatrix} \left.\begin{matrix} \overbrace{SD_1\ \ 0\ \cdots\ 0}^{|T|}\cdots\cdots\ 0 \\ 0\ \ SD_1\ \cdots\ 0\ \cdots\cdots 0 \\ \vdots \qquad \ddots \qquad \vdots \\ 0\cdots\cdots\cdots SD_1\cdots\cdots\ 0 \end{matrix}\right\}|T|(i=1) \\ \vdots \qquad \ddots \qquad \vdots \\ \left.\begin{matrix} 0\ \cdots\ \ SD_N\cdots\cdots\ 0 \\ 0\ \cdots\ \ 0\ SD_N\cdots 0 \\ \vdots \qquad \ddots \qquad \vdots \\ 0\cdots\cdots\ \ \underbrace{0\cdots\cdots\ SD_N}_{|T|} \end{matrix}\right\}|T|(i=N) \end{pmatrix}_{NT\times NT}$$

The following expressions are detailed descriptions of the abbreviations in Eq. (3d). *rampup* represents the coefficient of (1c) in the original model (1), represented by *rampup* in coefficient matrix $A$.

$$
rampup = \begin{pmatrix} \overbrace{\begin{matrix} 1 & \cdots & 0 \\ \vdots & \ddots & \vdots \\ 0 & \cdots & 1 \end{matrix}}^{|NT|} & \underbrace{\begin{matrix} 0 & \cdots & 0 \\ \vdots & \ddots & \vdots \\ 0 & \cdots & 0 \end{matrix}}_{N(2T-1)} \end{pmatrix}_{NT\times N(3T-1)} - (SR-G)\begin{pmatrix} \overbrace{0 \cdots 0}^{|NT|} & \overbrace{1\ 0 \cdots 0}^{|NT|} & \overbrace{0 \cdots 0}^{N(T-1)} \\ 0 \cdots 0 & 0\ 1 \cdots 0 & 0 \cdots 0 \\ \vdots & \ddots & \vdots \\ 0 \cdots 0 & 0 \cdots 0\ 1 & 0 \cdots 0 \end{pmatrix}_{NT\times N(3T-1)}
$$

$$
-\begin{pmatrix} \left.\begin{matrix} \overbrace{0\ 0 \cdots\ 0 \cdots 0\ \cdots\ 0}^{|NT|} & \overbrace{0\ \cdots\ 0}^{|N(2T-1)|} \\ 1\ 0 \cdots 0 \cdots 0\ \cdots\ 0 & 0\ \cdots\ 0 \\ \vdots\ \ \ddots\ \ \vdots & \vdots\ \ddots\ \vdots \\ \underbrace{0 \cdots 1}_{|T|}\ 0 \cdots 0\ \cdots\ 0 & 0\ \cdots\ 0 \end{matrix}\right\} T(i=1) \\ \vdots\ \ \ddots\ \ \vdots\ \ \vdots\ \ddots\ \vdots \\ \left.\begin{matrix} \overbrace{0\ \cdots\ 0 \cdots 0\ 0 \cdots\ 0}^{|NT|} & \overbrace{0\ \cdots\ 0}^{|N(2T-1)|} \\ 0\ \cdots\ 0 \cdots 1\ 0 \cdots 0 & 0\ \cdots\ 0 \\ \vdots\ \ \ddots\ \ \vdots & \vdots\ \ddots\ \vdots \\ 0\ \cdots\ 0 \cdots \underbrace{0 \cdots 1}_{|T|} & 0\ \cdots\ 0 \end{matrix}\right\} T(i=N) \end{pmatrix}_{NT\times N(3T-1)} - (UR-SR)\begin{pmatrix} \left.\begin{matrix} \overbrace{0 \cdots 0}^{|NT|} & \overbrace{0\ 0 \cdots 0 \cdots 0 \cdots 0}^{|NT|} & \overbrace{0 \cdots 0}^{N(T-1)} \\ 0 \cdots 0 & 1\ 0 \cdots 0 \cdots 0 \cdots 0 & 0 \cdots 0 \\ \vdots & \ddots & \vdots \\ 0 \cdots 0 & \underbrace{0 \cdots 1}_{|T|}\ 0 \cdots 0 \cdots 0 & 0 \cdots 0 \end{matrix}\right\} T(i=1) \\ \vdots\ \ \ddots\ \ \vdots \\ \left.\begin{matrix} \overbrace{0 \cdots 0}^{|NT|} & \overbrace{0 \cdots 0 \cdots 0\ 0 \cdots 0}^{|NT|} & \overbrace{0 \cdots 0}^{N(T-1)} \\ 0 \cdots 0 & 0 \cdots 0 \cdots 1\ 0 \cdots 0 & 0 \cdots 0 \\ \vdots & \ddots & \vdots \\ 0 \cdots 0 & 0 \cdots 0 \cdots \underbrace{0 \cdots 1}_{|T|} & 0 \cdots 0 \end{matrix}\right\} T(i=N) \end{pmatrix}_{NT\times N(3T-1)} \quad \text{(A.3)}
$$

In the following matrices, *shutdown* represents the coefficient of (1d) in the original model (1), represented by *shutdown* in coefficient matrix $A$. In particular, it has a row count of $N(T-1)$.

$$
shutdown = \begin{pmatrix} \left.\begin{matrix} \overbrace{1\ 0 \cdots 0 \cdots 0\ \cdots\ 0}^{|NT|} & \overbrace{0\ \cdots\ 0}^{|N(2T-1)|} \\ \vdots\ \ \ddots\ \ \vdots & \vdots\ \ddots\ \vdots \\ \underbrace{0 \cdots 1}_{|T|}\ 0 \cdots 0\ \cdots\ 0 & 0\ \cdots\ 0 \end{matrix}\right\} T-1(i=1) \\ \vdots\ \ \ddots\ \ \vdots\ \ \vdots\ \ddots\ \vdots \\ \left.\begin{matrix} \overbrace{0\ \cdots\ 0 \cdots 1\ 0 \cdots\ 0}^{|NT|} & \overbrace{0\ \cdots\ 0}^{|N(2T-1)|} \\ \vdots\ \ \ddots\ \ \vdots & \vdots\ \ddots\ \vdots \\ 0\ \cdots\ 0 \cdots \underbrace{0 \cdots 1}_{|T|} & 0\ \cdots\ 0 \end{matrix}\right\} T-1(i=N) \end{pmatrix}_{N(T-1)\times N(3T-1)} - SD'\begin{pmatrix} \left.\begin{matrix} \overbrace{0 \cdots 0}^{|NT|} & \overbrace{1\ 0 \cdots 0 \cdots 0 \cdots 0}^{|NT|} & \overbrace{0 \cdots 0}^{|N(T-1)|} \\ \vdots & \ddots & \vdots \\ 0 \cdots 0 & \underbrace{0 \cdots 1}_{|T|}\ 0 \cdots 0 \cdots 0 & 0 \cdots 0 \end{matrix}\right\} T-1(i=1) \\ \vdots\ \ \ddots\ \ \vdots \\ \left.\begin{matrix} \overbrace{0 \cdots 0}^{|NT|} & \overbrace{0 \cdots 0 \cdots 1\ 0 \cdots 0}^{|NT|} & \overbrace{0 \cdots 0}^{|N(T-1)|} \\ \vdots & \ddots & \vdots \\ 0 \cdots 0 & 0 \cdots 0 \cdots \underbrace{0 \cdots 1}_{|T|} & 0 \cdots 0 \end{matrix}\right\} T-1(i=N) \end{pmatrix}_{N(T-1)\times N(3T-1)}
$$

$$
-(G-SD)'\begin{pmatrix} \left.\begin{matrix} \overbrace{0 \cdots 0}^{|NT|} & \overbrace{0\ 1 \cdots 0 \cdots 0 \cdots 0}^{|NT|} & \overbrace{0 \cdots 0}^{|N(T-1)|} \\ \vdots & \ddots & \vdots \\ 0 \cdots 0 & \underbrace{0 \cdots 0\ 1}_{|T|} \cdots 0 \cdots 0 & 0 \cdots 0 \end{matrix}\right\} T-1(i=1) \\ \vdots\ \ \ddots\ \ \vdots \\ \left.\begin{matrix} \overbrace{0 \cdots 0}^{|NT|} & \overbrace{0 \cdots 0 \cdots 0\ 1 \cdots 0}^{|NT|} & \overbrace{0 \cdots 0}^{|N(T-1)|} \\ \vdots & \ddots & \vdots \\ 0 \cdots 0 & 0 \cdots 0 \cdots \underbrace{0 \cdots 0\ 1}_{|T|} & 0 \cdots 0 \end{matrix}\right\} T-1(i=N) \end{pmatrix}_{N(T-1)\times N(3T-1)} \quad \text{(A.4)}
$$

Similarly, in the following matrices, *rampdown* represents the coefficient of (1e) in the original model (1), represented by *rampdown* in coefficient matrix $A$.

$$rampdown=\begin{pmatrix} \left.\begin{matrix} \overbrace{0\,0\cdots\ 0\cdots0\ \cdots\ 0}^{|NT|}\ \overbrace{0\ \cdots\ 0}^{|N(2T-1)|} \\ 1\ 0\cdots0\cdots0\ \cdots\ 0\,0\ \cdots\ 0 \\ \vdots\quad\ddots\quad\vdots\,\vdots\ \ddots\ \vdots \\ \underbrace{0\cdots1}_{|T|}\ 0\cdots0\ \cdots\ 0\,0\ \cdots\ 0 \end{matrix}\right\}T(i=1) \\ \vdots\quad\ddots\quad\vdots\,\vdots\ \ddots\ \vdots \\ \left.\begin{matrix} \overbrace{0\ \cdots\ 0\cdots0\ 0\cdots\ 0}^{|NT|}\ \overbrace{0\ \cdots\ 0}^{|N(2T-1)|} \\ 0\ \cdots\ 0\cdots1\ 0\cdots0\,0\ \cdots\ 0 \\ \vdots\quad\ddots\quad\vdots\,\vdots\ \ddots\ \vdots \\ 0\ \cdots\ 0\cdots\underbrace{0\cdots1}_{|T|}\ 0\ 0\ \cdots\ 0 \end{matrix}\right\}T(i=N) \end{pmatrix}_{NT\times N(3T-1)} -(SD-G)\begin{pmatrix} \left.\begin{matrix} \overbrace{0\ \cdots\ 0}^{|NT|}\ \overbrace{0\ 0\cdots0\cdots0\ \cdots0}^{|NT|}\ \overbrace{0\ \cdots\ 0}^{N(T-1)} \\ 0\ \cdots\ 0\,1\ 0\cdots0\cdots0\ \cdots\ 0\,0\ \cdots\ 0 \\ \vdots\quad\ddots\quad\vdots \\ 0\cdots\ 0\,\underbrace{0\ \cdots1}_{|T|}\,0\cdots0\ \cdots0\,0\ \cdots\ 0 \end{matrix}\right\}T(i=1) \\ \vdots\quad\ddots\quad\vdots \\ \left.\begin{matrix} \overbrace{0\ \cdots\ 0}^{|NT|}\ \overbrace{0\cdots0\cdots0\,0\cdots0}^{|NT|}\ \overbrace{0\ \cdots\ 0}^{N(T-1)} \\ 0\ \cdots\ 0\,0\cdots0\cdots1\ 0\cdots\ 0\,0\ \cdots\ 0 \\ \vdots\quad\ddots\quad\vdots \\ 0\cdots\ 0\,0\cdots0\cdots\underbrace{0\ \cdots1}_{|T|}\ 0\ 0\ \cdots\ 0 \end{matrix}\right\}T(i=N) \end{pmatrix}_{NT\times N(3T-1)} \tag{A.5}$$

$$-(DR-SD)\begin{pmatrix} \overbrace{0\ \cdots\ 0}^{|NT|}\ \overbrace{1\ 0\cdots0}^{|NT|}\ \overbrace{0\ \cdots\ 0}^{N(T-1)} \\ 0\ \cdots\ 0\,0\ 1\cdots0\,0\ \cdots\ 0 \\ \vdots\quad\ddots\quad\vdots \\ 0\cdots\ 0\,0\ \cdots0\,1\,0\ \cdots\ 0 \end{pmatrix}_{NT\times N(3T-1)} -\begin{pmatrix} \overbrace{1\ \ldots\ 0}^{|NT|}0\ \ldots\ 0 \\ \vdots\ \ddots\ \vdots\,\vdots\ \ddots\ \vdots \\ 0\ \cdots\ 1\underbrace{0\ \cdots\ 0}_{N(2T-1)} \end{pmatrix}_{NT\times N(3T-1)}$$

Ramp-up, shutdown, and ramp-down are represented as sub-parts of matrix A by (1c) -(1e). These are obtained from the variables coefficient matrix on the left side of the original constraint, which shows that variables are controlled by startup ramp rates and ramp-up rates in (1c), as well as shutdown ramp rates and ramp-down rates in (1d) and (1e).

The following *state transition* matrix represents the coefficient of (1f) in the original model (1), represented by *state transition* in coefficient matrix $A$.

$$state\ transition=\begin{pmatrix} \left.\begin{matrix} 0\ \cdots\cdots\ 0\ \overbrace{-1\ 1\ 0\cdots\cdots0}^{|T|}\ 0\cdots\cdots\ 0\ \overbrace{-1\ 0\ \cdots\cdots0}^{|T-1|}\ 0\cdots\cdots\ 0 \\ 0\ \cdots\cdots\ 0\,0\ -1\,1\,0\cdots0\ 0\cdots\cdots\ 0\ 0\ -1\ 0\cdots0\ 0\cdots\cdots\ 0 \\ \vdots\qquad\ddots\qquad\vdots \\ \underbrace{0\ \cdots\cdots\ 0}_{|NT|}\ \underbrace{0\ \cdots\ 0\ -1\ 1\ 0\cdots\cdots\ 0}_{|NT|}\ \underbrace{0\ \cdots\ 0\ \ -1\,0\ \cdots\cdots\ 0}_{|N(T-1)|} \end{matrix}\right\}t=2,\cdots T\ (i=1) \\ \vdots\qquad\ddots\qquad\vdots \\ \left.\begin{matrix} 0\ \cdots\cdots\ 0\,0\ \cdots\cdots\ 0\ \overbrace{-1\ 1\ 0\cdots\cdots0}^{|T|}\ 0\ \cdots\cdots\ 0\ \overbrace{-1\ 0\cdots\cdots0}^{|T-1|} \\ 0\ \cdots\cdots\ 0\,0\ \cdots\cdots\ 0\,0\ -1\,1\,0\cdots0\,0\ \cdots\cdots\ 0\,0\ \ -1\,0\cdots0 \\ \vdots\qquad\ddots\qquad\vdots \\ \underbrace{0\ \cdots\cdots\ 0}_{|NT|}\ \underbrace{0\ \cdots\cdots\ 0\,0\cdots0\ -1\ 1}_{|NT|}\ \underbrace{0\ \cdots\cdots0\,0\cdots0\ \ -1}_{|N(T-1)|} \end{matrix}\right\}t=2,\cdots T\ (i=N) \end{pmatrix}_{N(T-1)\times N(3T-1)} \tag{A.6}$$

Next, both $T^{on}$ and $T^{off}$ matrices are parameters that are unaffected by the optimization periods for any problems that are constrained by physical constraints of generators.

$$T^{on}=\begin{pmatrix} T_1^{on} & \ldots & 0 \\ \vdots & \ddots & \vdots \\ 0 & \cdots & T_N^{on} \end{pmatrix}_{NT\times NT} \quad T^{off}=\begin{pmatrix} T_1^{off} & \ldots & 0 \\ \vdots & \ddots & \vdots \\ 0 & \cdots & T_N^{off} \end{pmatrix}_{NT\times NT} \tag{A.7}$$

The following $\min on$ matrix represents the coefficient of (1g) in the original model (1), represented by $\min on$ in coefficient matrix $A$.

$$\min on = \left(\begin{array}{l} \left.\begin{array}{l} 0\ \cdots\cdots 0\ \overbrace{0\cdots -1_{T_1^{on}+1}\cdots\cdots 0}^{|T|}\ \cdots\cdots 0\ \overbrace{1\ \cdots 1_{T_1^{on}+1}\ \cdots\cdots\ 0}^{|T-1|}\cdots\cdots\ 0 \\ 0\ \cdots\cdots 0\ 0\ \cdots\ -1_{T_1^{on}+2}\cdots\ 0\ \cdots\cdots 0\ 0\ \ 1\cdots 1_{T_1^{on}+2}\cdots\cdots 0\cdots\cdots\ 0 \\ \quad\vdots \qquad\qquad \ddots \qquad\qquad \vdots \\ 0\ \cdots\cdots 0\ \underbrace{0\ \cdots\ 0\cdots\cdots -1_T\ \ 0\ \cdots\cdots 0}_{|NT|}\ \underbrace{0\cdots\cdots 1_{T-T_1^{on}+1}\cdots 1_T\cdots\cdots 0}_{|N(T-1)|} \end{array}\right\} t = T_1^{on}+1,\cdots T\ \ (i=1) \\ \quad\vdots \qquad\qquad\qquad \ddots \qquad\qquad\qquad \vdots \\ \left.\begin{array}{l} 0\ \cdots\cdots 0\ 0\cdots\cdots \overbrace{0\cdots -1_{T_N^{on}+1}\cdots\cdots 0}^{|T|}\ 0\cdots\cdots 0\ \overbrace{1\ \cdots 1_{T_N^{on}+1}\ \cdots\cdots\ 0}^{|T-1|} \\ 0\ \cdots\cdots 0\ 0\cdots\cdots 0\ \cdots -1_{T_N^{on}+2}\cdots 0\ \ 0\cdots\cdots 0\ 0\ \ 1\cdots 1_{T_N^{on}+2}\cdots\cdots 0 \\ \quad\vdots \qquad\qquad \ddots \qquad\qquad \vdots \\ \underbrace{0\ \cdots\cdots 0}_{|NT|}\ \underbrace{0\cdots\cdots 0\ \cdots\ 0\cdots\cdots -1_T}_{|NT|}\ \underbrace{0\cdots\cdots 0\ 0\cdots\cdots 1_{T-T_N^{on}+1}\cdots 1_T}_{|N(T-1)|} \end{array}\right\} t = T_N^{on}+1,\cdots T\ \ (i=N) \end{array}\right)_{(NT-\sum_{i=1}^{N} T_i^{on})\times N(3T-1)} \tag{A.8}$$

The matrix $\min off$ represents the coefficient of (1h) in the original model (1), represented by $\min off$ in coefficient matrix $A$.

$$\min off = \left(\begin{array}{l} \left.\begin{array}{l} 0\ \cdots\cdots 0\ \overbrace{1\cdots\cdots\cdots\cdots\ 0}^{|T|}\cdots\cdots 0\ \overbrace{1\ \cdots 1_{T_1^{off}+1}\ \cdots\cdots\ 0}^{|T-1|}\cdots\cdots 0 \\ 0\ \cdots\cdots 0\ 0\ 1\ \ \cdots\cdots\ \ 0\cdots\cdots 0\ \ 0\ 1\cdots 1_{T_1^{off}+2}\cdots\cdots 0\cdots\cdots 0 \\ \quad\vdots \qquad\qquad \ddots \qquad\qquad \vdots \\ 0\ \cdots\cdots 0\ \underbrace{0\cdots 1_{T-T_1^{off}}\cdots 0\cdots\cdots 0}_{|NT|}\ \underbrace{0\cdots\cdots 1_{T-T_1^{off}+1}\cdots 1_T\cdots\cdots 0}_{|N(T-1)|} \end{array}\right\} t = T_1^{off}+1,\cdots T\ \ (i=1) \\ \quad\vdots \qquad\qquad\qquad \ddots \qquad\qquad\qquad \vdots \\ \left.\begin{array}{l} 0\ \cdots\cdots 0\ 0\cdots\cdots \overbrace{1\ \cdots\cdots\cdots\cdots\ 0}^{|T|}\ 0\cdots\cdots \overbrace{1\ \cdots 1_{T_N^{off}+1}\ \cdots\cdots\ 0}^{|T-1|} \\ 0\ \cdots\cdots 0\ 0\cdots\cdots 0\ 1\cdots\cdots\cdots 0\ \ 0\ 0\cdots\cdots 0\ \ 1\cdots 1_{T_N^{off}+2}\cdots\cdots 0 \\ \quad\vdots \qquad\qquad \ddots \qquad\qquad \vdots \\ \underbrace{0\ \cdots\cdots 0}_{|NT|}\ \underbrace{0\cdots\cdots 0\ \cdots 1_{T-T_1^{off}}\cdots 0}_{|NT|}\ \underbrace{0\cdots\cdots 0\cdots\cdots 1_{T-T_N^{off}+1}\cdots 1_T}_{|N(T-1)|} \end{array}\right\} t = T_N^{off}+1,\cdots T\ \ (i=N) \end{array}\right)_{(NT-\sum_{i=1}^{N} T_i^{off})\times N(3T-1)} \tag{A.9}$$

State transition, min on, and min off are constraints from (1f) to (1h) expressed separately as sub-parts of matrix A. They are computed by using the variables coefficient matrix on the left side of the original constraint, which represents the minimum up (on) and down (off) time limits.

**Appendix B: Three periods Case Study with Synthetic Data**

Set startup ramp limit, ramp-up limit, shutdown limit, ramp-down limit and maximum power generation capacity are 20 at every generator in three periods. Minimum on (up) time and minimum off (down) time are 2 for two generators, $G_{i0}$ and $U_{i0}$ are 0. The cost coefficients in the objective function are 3,7,5,4,6,4,1,1,1,1,1,1,18,16,14,20.

Generation shift factor matrix is $\begin{pmatrix} 0.8 & 0 & 0 & 0.7 & 0 & 0 \\ 0 & 0.8 & 0 & 0 & 0.7 & 0 \\ 0 & 0 & 0.8 & 0 & 0 & 0.7 \\ 0.6 & 0 & 0 & 0.9 & 0 & 0 \\ 0 & 0.6 & 0 & 0 & 0.9 & 0 \\ 0 & 0 & 0.6 & 0 & 0 & 0.9 \end{pmatrix}$,demand matrix is $D = \begin{pmatrix} 8 \\ 7 \\ 12 \\ 10 \\ 8 \\ 8 \end{pmatrix}$.

To perform this problem via parametric analysis for the MILP approach, we first reformat the three-period case, including two generators, two lines, and two buses. Then, we adjust the capacity of the first line, assigning the uncertainty parameter theta, which varies between 0 and 10.

$$z(\theta) = \min_{x,y} 3G_{11} + 7G_{12} + 5G_{13} + 4G_{21} + 6G_{22} + 4G_{23} + 1U_{11} + 1U_{12} + 1U_{13}$$
$$+ 1U_{21} + 1U_{22} + 1U_{23} + 18V_{12} + 16V_{13} + 14V_{22} + 20V_{23}$$

$$s.t. \begin{cases} G_{11} + G_{21} = 18 \\ G_{12} + G_{22} = 15 \\ G_{13} + G_{23} = 20 \\ 0.8G_{11} + 0.7G_{21} \le 13.2 + \theta_1 \\ 0.8G_{12} + 0.7G_{22} \le 11 + \theta_1 \\ 0.8G_{13} + 0.7G_{23} \le 15 + \theta_1 \\ 0.6G_{11} + 0.9G_{21} \le 14.4 \\ 0.6G_{12} + 0.9G_{22} \le 12 \\ 0.6G_{13} + 0.9G_{23} \le 15 \\ G_{11} \le 20 \\ G_{12} - G_{11} \le 20 \\ G_{13} - G_{12} \le 20 \end{cases} \begin{cases} G_{21} \le 20 \\ G_{22} - G_{21} \le 20 \\ G_{23} - G_{22} \le 20 \\ G_{11} - 20U_{11} \le 0 \\ G_{12} - 20U_{12} \le 0 \\ G_{21} - 20U_{21} \le 0 \\ G_{22} - 20U_{22} \le 0 \\ -G_{11} \le 20 \\ G_{11} - G_{12} \le 20 \\ G_{12} - G_{13} \le 20 \\ -G_{21} \le 20 \\ G_{21} - G_{22} \le 20 \\ G_{22} - G_{23} \le 20 \end{cases} \begin{cases} U_{12} - U_{11} - V_{12} \le 0 \\ U_{22} - U_{21} - V_{22} \le 0 \\ U_{13} - U_{12} - V_{13} \le 0 \\ U_{23} - U_{22} - V_{23} \le 0 \\ V_{12} + V_{13} - U_{13} \le 0 \\ V_{22} + V_{23} - U_{23} \le 0 \\ V_{12} + V_{13} + U_{11} \le 1 \\ V_{22} + V_{23} + U_{21} \le 1 \\ U_{it} \in \{0,1\} \\ V_{it} \in \{0,1\} \\ 0 \le \theta_s \le 10,\ s = 1 \\ 0 \le G_{it} \le 20 \end{cases}$$

**Step0:** Initialize current upper bound $\bar{z}_{upper} = \infty$.

**Step1**: Solve the fully relaxed problem as follows at root node (i.e., node 0).

(1) When $\theta_1 = 0$,we can get following optimal solutions:

$$\begin{cases} G_{11}=6\ G_{12}=5\ G_{13}=10\ G_{21}=12\ \ G_{22}=10\ G_{23}=10 \\ U_{11}=0.5\ U_{12}=0.5\ U_{13}=0.5\ U_{21}=0.6\ \ U_{22}=0.5\ U_{23}=0.5 \\ V_{12}=0\ V_{13}=0\ V_{22}=0\ V_{23}=0 \end{cases} \tag{B.1}$$

The corresponding optimal Lagrange multipliers are shown as follows:

$$\begin{cases} \mu_1=10.5\ \mu_5=10/3\ \mu_6=3.5\ \mu_{23}=1\ \mu_{24}=2\ \mu_{26}=1 \\ \mu_{39}=0.15\ \mu_{40}=0.05\ \mu_{42}=0.1\ \mu_{55}=17\ \mu_{56}=14\ \mu_{57}=14\ \mu_{58}=19 \\ \lambda_1=-11.4\ \lambda_2=-9\ \lambda_3=-7.25 \end{cases} \tag{B.2}$$

Recall Section 4 and Section 5.1; using the results in Eq. (B2), we have the following basic matrix derived from KKT conditions for this problem.

$$A_p=\begin{bmatrix} 0.8 & 0 & 0 & 0.7 & 0 & 0 & 0 & 0 & 0 & 0 & 0 & 0 & 0 & 0 & 0 & 0 \\ 0 & 0.6 & 0 & 0 & 0.9 & 0 & 0 & 0 & 0 & 0 & 0 & 0 & 0 & 0 & 0 & 0 \\ 0 & 0 & 0.6 & 0 & 0 & 0.9 & 0 & 0 & 0 & 0 & 0 & 0 & 0 & 0 & 0 & 0 \\ 0 & 0 & 0 & 0 & 0 & 0 & -1 & 1 & 0 & 0 & 0 & 0 & -1 & 0 & 0 & 0 \\ 0 & 0 & 0 & 0 & 0 & 0 & 0 & -1 & 1 & 0 & 0 & 0 & 0 & -1 & 0 & 0 \\ 0 & 0 & 0 & 0 & 0 & 0 & 0 & 0 & 0 & -1 & 1 & 0 & 0 & 0 & -1 \\ 0 & 0 & 1 & 0 & 0 & 0 & 0 & 0 & -20 & 0 & 0 & 0 & 0 & 0 & 0 & 0 \\ 0 & 0 & 0 & 1 & 0 & 0 & 0 & 0 & 0 & -20 & 0 & 0 & 0 & 0 & 0 & 0 \\ 0 & 0 & 0 & 0 & 0 & 1 & 0 & 0 & 0 & 0 & 0 & -20 & 0 & 0 & 0 & 0 \\ 0 & 0 & 0 & 0 & 0 & 0 & 0 & 0 & 0 & 0 & 0 & 0 & -1 & 0 & 0 & 0 \\ 0 & 0 & 0 & 0 & 0 & 0 & 0 & 0 & 0 & 0 & 0 & 0 & 0 & -1 & 0 & 0 \\ 0 & 0 & 0 & 0 & 0 & 0 & 0 & 0 & 0 & 0 & 0 & 0 & 0 & 0 & -1 & 0 \\ 0 & 0 & 0 & 0 & 0 & 0 & 0 & 0 & 0 & 0 & 0 & 0 & 0 & 0 & 0 & -1 \\ 1 & 0 & 0 & 1 & 0 & 0 & 0 & 0 & 0 & 0 & 0 & 0 & 0 & 0 & 0 & 0 \\ 0 & 1 & 0 & 0 & 1 & 0 & 0 & 0 & 0 & 0 & 0 & 0 & 0 & 0 & 0 & 0 \\ 0 & 0 & 1 & 0 & 0 & 1 & 0 & 0 & 0 & 0 & 0 & 0 & 0 & 0 & 0 & 0 \end{bmatrix} \tag{B.3}$$

Similar to Section 5.1, the optimal decision variable $\omega_{0,1}^*$ can be obtained from Eq. (15), the critical region $CR_{0,1}$ obtained from Eq. (24), and the optimal value $z_{0,1}$ derived from Eq. (16).

$$\omega_{0,1}^*=A_p^{-1}(b_p+F_p^{\max}+\theta_p)=\left(10\cdot\theta_1+6,5,10,12-10\cdot\theta_1,10,10,0.5,0.5,0.5,0.6-0.5\cdot\theta_1,0.5,0.5,0,0,0,0\right)^T \tag{B.4}$$

$$\begin{aligned} &CR_{0,1}=\left\{\theta\in\Theta\,|\,A_s\omega^*\le b_s+\hat{F}_s^{\max}\right\}=\left\{0\le\theta_1\le0.2\right\} \\ &z_{0,1}=254.1-10.5\cdot\theta_1 \end{aligned} \tag{B.5}$$

(2) When $\theta_1=0.2$, the twenty-sixth inactive constraint in all constraints of the problem, changes to active constraint; the fourteenth active constraint changes to inactive constraint, the basic matrix $A_p$ will change with it. When $\theta_1=0.2$,we can get following optimal solutions:

$$\begin{cases} G_{11}=8\ G_{12}=5\ G_{13}=10\ G_{21}=10\ \ G_{22}=10\ G_{23}=10 \\ U_{11}=0.5\ U_{12}=0.5\ U_{13}=0.5\ U_{21}=0.5\ \ U_{22}=0.5\ U_{23}=0.5 \\ V_{12}=0\ V_{13}=0\ V_{22}=0\ V_{23}=0 \end{cases} \tag{B.6}$$

The associated optimal Lagrange multipliers are:

$$\begin{cases} \mu_1=10\ \mu_5=10/3\ \mu_6=10/3\ \mu_{23}=1\ \mu_{24}=2\ \mu_{25}=1 \\ \mu_{26}=2\ \mu_{39}=0.15\ \mu_{42}=0.15\ \mu_{55}=17\ \mu_{56}=14\ \mu_{57}=13\ \mu_{58}=18 \\ \lambda_1=-11\ \lambda_2=-9\ \ \lambda_3=-7.15 \end{cases} \tag{B.7}$$

When $\theta_1 = 0.2$ , we also get the following updated basic matrix.

$$A_p = \begin{bmatrix} 0.8 & 0 & 0 & 0.7 & 0 & 0 & 0 & 0 & 0 & 0 & 0 & 0 & 0 & 0 & 0 & 0 \\ 0 & 0.6 & 0 & 0 & 0.9 & 0 & 0 & 0 & 0 & 0 & 0 & 0 & 0 & 0 & 0 & 0 \\ 0 & 0 & 0.6 & 0 & 0 & 0.9 & 0 & 0 & 0 & 0 & 0 & 0 & 0 & 0 & 0 & 0 \\ 0 & 0 & 0 & 0 & 0 & 0 & -1 & 1 & 0 & 0 & 0 & 0 & -1 & 0 & 0 & 0 \\ 0 & 0 & 0 & 0 & 0 & 0 & 0 & -1 & 1 & 0 & 0 & 0 & 0 & -1 & 0 & 0 \\ 0 & 0 & 0 & 0 & 0 & 0 & 0 & 0 & -1 & 1 & 0 & 0 & 0 & -1 & 0 \\ 0 & 0 & 0 & 0 & 0 & 0 & 0 & 0 & 0 & -1 & 1 & 0 & 0 & 0 & -1 \\ 0 & 0 & 1 & 0 & 0 & 0 & 0 & 0 & -20 & 0 & 0 & 0 & 0 & 0 & 0 & 0 \\ 0 & 0 & 0 & 0 & 0 & 1 & 0 & 0 & 0 & 0 & 0 & -20 & 0 & 0 & 0 & 0 \\ 0 & 0 & 0 & 0 & 0 & 0 & 0 & 0 & 0 & 0 & 0 & -1 & 0 & 0 & 0 \\ 0 & 0 & 0 & 0 & 0 & 0 & 0 & 0 & 0 & 0 & 0 & 0 & -1 & 0 & 0 \\ 0 & 0 & 0 & 0 & 0 & 0 & 0 & 0 & 0 & 0 & 0 & 0 & 0 & -1 & 0 \\ 0 & 0 & 0 & 0 & 0 & 0 & 0 & 0 & 0 & 0 & 0 & 0 & 0 & 0 & -1 \\ 1 & 0 & 0 & 1 & 0 & 0 & 0 & 0 & 0 & 0 & 0 & 0 & 0 & 0 & 0 & 0 \\ 0 & 1 & 0 & 0 & 1 & 0 & 0 & 0 & 0 & 0 & 0 & 0 & 0 & 0 & 0 & 0 \\ 0 & 0 & 1 & 0 & 0 & 1 & 0 & 0 & 0 & 0 & 0 & 0 & 0 & 0 & 0 & 0 \end{bmatrix} \tag{B.8}$$

Similarly to the first region computing process, the optimal decision variable $\omega_{0,2}^*$ can be obtained from Eq. (15), the optimal value $z_{0,2}$ derived via Eq. (16), and the critical region $CR_{0,2}$ obtained using Eq. (24).

$$\omega_{0,2}^* = A_p^{-1}(b_p + F_p^{\max} + \theta_p) = (\ 10\cdot\theta_1+6,5,10,12-10\cdot\theta_1,10,10,0.5,0.5,0.5,0.5,0.5,0.5,0,0,0,0)^T \tag{B.9}$$

$$\begin{aligned} &CR_{0,2} = \left\{\theta\in\Theta \mid A_s\omega^* \le b_s + \hat{F}_s^{\max}\ \right\} = \left\{0.2 \le \theta_1 \le 0.4\right\} \\ &z_{0,2}=254-10*\theta_1 \end{aligned} \tag{B.10}$$

(3) When $\theta_1 = 0.4$ , the thirteenth inactive constraint in this problem, changes to active constraint, the twenty-fifth active constraint changes to inactive constraint, the part of $A_p$ will change with it.

When $\theta_1 = 0.4$ ,we can get following optimal solutions:

$$\begin{cases} G_{11}=10\ G_{12}=5\ G_{13}=10\ G_{21}=8\ \ G_{22}=10\ G_{23}=10 \\ U_{11}=0.5\ U_{12}=0.5\ U_{13}=0.5\ U_{21}=0.5\ \ U_{22}=0.5\ U_{23}=0.5 \\ V_{12}=0\ V_{13}=0\ V_{22}=0\ V_{23}=0 \end{cases} \tag{B.11}$$

The corresponding optimal Lagrange multipliers are

$$\begin{cases} \mu_1=9.5\ \mu_5=10/3\ \mu_6=19/6\ \ \mu_{24}=1\ \ \mu_{25}=1\ \mu_{26}=2 \\ \mu_{37}=0.05\ \mu_{39}=0.1\ \mu_{42}=0.15\ \mu_{55}=18\ \mu_{56}=15\ \ \mu_{57}=13\ \mu_{58}=18 \\ \lambda_1=-10.65\ \lambda_2=-9\ \ \lambda_3=-7 \end{cases} \tag{B.12}$$

When $\theta_1=0.4$, we also get the following updated basic matrix.

$$A_p=\begin{bmatrix} 0.8 & 0 & 0 & 0.7 & 0 & 0 & 0 & 0 & 0 & 0 & 0 & 0 & 0 & 0 & 0 & 0 \\ 0 & 0.6 & 0 & 0 & 0.9 & 0 & 0 & 0 & 0 & 0 & 0 & 0 & 0 & 0 & 0 & 0 \\ 0 & 0 & 0.6 & 0 & 0 & 0.9 & 0 & 0 & 0 & 0 & 0 & 0 & 0 & 0 & 0 & 0 \\ 0 & 0 & 0 & 0 & 0 & 0 & 0 & -1 & 1 & 0 & 0 & 0 & 0 & -1 & 0 & 0 \\ 0 & 0 & 0 & 0 & 0 & 0 & 0 & 0 & 0 & -1 & 1 & 0 & 0 & 0 & -1 & 0 \\ 0 & 0 & 0 & 0 & 0 & 0 & 0 & 0 & 0 & 0 & -1 & 1 & 0 & 0 & 0 & -1 \\ 1 & 0 & 0 & 0 & 0 & 0 & -20 & 0 & 0 & 0 & 0 & 0 & 0 & 0 & 0 & 0 \\ 0 & 0 & 1 & 0 & 0 & 0 & 0 & 0 & -20 & 0 & 0 & 0 & 0 & 0 & 0 & 0 \\ 0 & 0 & 0 & 0 & 0 & 1 & 0 & 0 & 0 & 0 & 0 & -20 & 0 & 0 & 0 & 0 \\ 0 & 0 & 0 & 0 & 0 & 0 & 0 & 0 & 0 & 0 & 0 & 0 & -1 & 0 & 0 & 0 \\ 0 & 0 & 0 & 0 & 0 & 0 & 0 & 0 & 0 & 0 & 0 & 0 & 0 & -1 & 0 & 0 \\ 0 & 0 & 0 & 0 & 0 & 0 & 0 & 0 & 0 & 0 & 0 & 0 & 0 & 0 & -1 & 0 \\ 0 & 0 & 0 & 0 & 0 & 0 & 0 & 0 & 0 & 0 & 0 & 0 & 0 & 0 & 0 & -1 \\ 1 & 0 & 0 & 1 & 0 & 0 & 0 & 0 & 0 & 0 & 0 & 0 & 0 & 0 & 0 & 0 \\ 0 & 1 & 0 & 0 & 1 & 0 & 0 & 0 & 0 & 0 & 0 & 0 & 0 & 0 & 0 & 0 \\ 0 & 0 & 1 & 0 & 0 & 1 & 0 & 0 & 0 & 0 & 0 & 0 & 0 & 0 & 0 & 0 \end{bmatrix} \tag{B.13}$$

Similarly, the optimal decision variable $\omega_{0,3}^*$ can be obtained from Eq. (15), the optimal value $z_{0,3}$ derived from Eq. (16), and the critical region $CR_{0,3}$ derived from Eq. (24).

$$\begin{cases} \omega_{0,3}^* = A_p^{-1}(b_p + F_p^{\max} + \theta_p) = (\ 10\cdot\theta_1+6,5,10,12-10\cdot\theta_1,10,10,0.5\cdot\theta_1+0.3,0.5,0.5,0.5,0.5,0.5,0,0,0,0)^T \\ CR_{0,3} = \left\{\theta\in\Theta \mid A_s\omega^* \le b_s + \hat{F}_s^{\max}\ \right\} = \{0.4 \le \theta_1 \le 1.2\} \\ z_{0,3}=253.8-9.5\cdot\theta_1 \end{cases} \tag{B.14}$$

(4) When $\theta_1=1.2$, the thirty-seventh inactive constraint in this problem, changes to active constraint, the first active constraint in all constraints changes to inactive constraint, the part of $A_p$ will change with it.

Meanwhile, we can get following optimal solutions:

$$\begin{cases} G_{11}=18\ G_{12}=5\ G_{13}=10\ G_{21}=0\ \ G_{22}=10\ G_{23}=10 \\ U_{11}=0.9\ U_{12}=0.5\ U_{13}=0.5\ U_{21}=0.5\ \ U_{22}=0.5\ U_{23}=0.5 \\ V_{12}=0\ V_{13}=0\ V_{22}=0\ V_{23}=0 \end{cases} \tag{B.15}$$

The optimal Lagrange multipliers are

$$\begin{cases} \mu_5=10/3\ \mu_6=19/6\ \mu_{24}=1\ \mu_{25}=1\ \mu_{26}=2\ \mu_{34}=0.95 \\ \mu_{37}=0.05\ \mu_{39}=0.1\ \mu_{42}=0.15\ \mu_{55}=18\ \mu_{56}=15\ \mu_{57}=13\ \mu_{58}=18 \\ \lambda_1=-3.05\ \lambda_2=-9\ \ \lambda_3=-7 \end{cases} \tag{B16}$$

At this time, we have the following basic matrix derived from KKT conditions.

$$A_p=\begin{bmatrix} 0 & 0.6 & 0 & 0 & 0.9 & 0 & 0 & 0 & 0 & 0 & 0 & 0 & 0 & 0 & 0 & 0 \\ 0 & 0 & 0.6 & 0 & 0 & 0.9 & 0 & 0 & 0 & 0 & 0 & 0 & 0 & 0 & 0 & 0 \\ 0 & 0 & 0 & 0 & 0 & 0 & 0 & -1 & 1 & 0 & 0 & 0 & 0 & -1 & 0 & 0 \\ 0 & 0 & 0 & 0 & 0 & 0 & 0 & 0 & 0 & -1 & 1 & 0 & 0 & 0 & -1 & 0 \\ 0 & 0 & 0 & 0 & 0 & 0 & 0 & 0 & 0 & 0 & -1 & 1 & 0 & 0 & 0 & -1 \\ 0 & 0 & 0 & -1 & 0 & 0 & 0 & 0 & 0 & 0 & 0 & 0 & 0 & 0 & 0 & 0 \\ 1 & 0 & 0 & 0 & 0 & 0 & -20 & 0 & 0 & 0 & 0 & 0 & 0 & 0 & 0 & 0 \\ 0 & 0 & 1 & 0 & 0 & 0 & 0 & 0 & -20 & 0 & 0 & 0 & 0 & 0 & 0 & 0 \\ 0 & 0 & 0 & 0 & 0 & 1 & 0 & 0 & 0 & 0 & 0 & -20 & 0 & 0 & 0 & 0 \\ 0 & 0 & 0 & 0 & 0 & 0 & 0 & 0 & 0 & 0 & 0 & 0 & -1 & 0 & 0 & 0 \\ 0 & 0 & 0 & 0 & 0 & 0 & 0 & 0 & 0 & 0 & 0 & 0 & 0 & -1 & 0 & 0 \\ 0 & 0 & 0 & 0 & 0 & 0 & 0 & 0 & 0 & 0 & 0 & 0 & 0 & 0 & -1 & 0 \\ 0 & 0 & 0 & 0 & 0 & 0 & 0 & 0 & 0 & 0 & 0 & 0 & 0 & 0 & 0 & -1 \\ 1 & 0 & 0 & 1 & 0 & 0 & 0 & 0 & 0 & 0 & 0 & 0 & 0 & 0 & 0 & 0 \\ 0 & 1 & 0 & 0 & 1 & 0 & 0 & 0 & 0 & 0 & 0 & 0 & 0 & 0 & 0 & 0 \\ 0 & 0 & 1 & 0 & 0 & 1 & 0 & 0 & 0 & 0 & 0 & 0 & 0 & 0 & 0 & 0 \end{bmatrix} \tag{B.17}$$

The optimal decision variable $\omega_{0,4}^{*}$ can be obtained from Eq. (15), the optimal value $z_{0,4}$ derived from Eq. (16), and the critical region $CR_{0,4}$ derived from Eq. (24).

At the root node, the entire parameter space is now covered.

$$\begin{cases} \omega_{0,4}^{*}=A_p^{-1}(b_p+F_p^{\max}+\theta_p)=\left(18\ \ 5\ \ 10\ \ 0\ \ 10\ \ 10\ \ 0.9\ \ 0.5\ \ 0.5\ \ 0.5\ \ 0.5\ \ 0.5\ \ 0\ \ 0\ \ 0\ \ 0\right)^{T} \\ CR_{0,4}=\left\{\theta\in\Theta \mid A_s\omega^{*}\le b_s+\hat{F}_s^{\max}\right\}=\left\{1.2\le\theta_1\le 10\right\} \\ z_{0,4}=242.4 \end{cases} \tag{B.18}$$

The detailed process of the root node full relaxation problem is described above. Then, because the method is very similar and the details are not presented here, we should calculate some other feasible nodes. The optimal values determined for each node, as well as the details of their critical regions, are shown in Table 1. We set the $\delta=0$, the accuracy rate is 100%.

**Table 1: Details Concerning Potential Node Parametric Analysis**

| Node 1. $U_{11}=1$ | Node 2. $U_{11}=0$ |
|---|---|
| $\begin{cases} z=254.6-10.5\cdot\theta_1, & 0\le\theta_1\le 0.2 \\ z=254.5-10\cdot\theta_1, & 0.2\le\theta_1\le 1.2 \\ z=242.5, & 1.2\le\theta_1\le 10 \end{cases}$ | Infeasible |
| Node 3. $U_{11}=1\ U_{21}=1$ | $\begin{cases} z=255-10\cdot\theta_1, & 0\le\theta_1\le 1.2 \\ z=243, & 1.2\le\theta_1\le 10 \end{cases}$ |
| Node 4. $U_{11}=1\ U_{21}=0$ | Infeasible |
| Node 5. $U_{11}=1\ U_{21}=1\ U_{12}=1$ | $\begin{cases} z=255.5-10\cdot\theta_1, & 0\le\theta_1\le 1.2 \\ z=243.5, & 1.2\le\theta_1\le 10 \end{cases}$ |
| Node 6. $U_{11}=1\ U_{21}=1\ U_{12}=0$ | Infeasible |
| Node 7. $U_{11}=1\ U_{21}=1\ U_{12}=1\ U_{22}=1$ | $\begin{cases} z=256-10\cdot\theta_1, & 0\le\theta_1\le 1.2 \\ z=244, & 1.2\le\theta_1\le 10 \end{cases}$ |
| Node 8. $U_{11}=1\ U_{21}=1\ U_{12}=1\ U_{22}=0$ | Infeasible |
| Node 9. $U_{11}=1\ U_{21}=1\ U_{12}=1\ U_{22}=1\ U_{13}=1$ | $\begin{cases} z=256.5-10\cdot\theta_1, & 0\le\theta_1\le 1.2 \\ z=244.5, & 1.2\le\theta_1\le 10 \end{cases}$ |
| Node 10. $U_{11}=1\ U_{21}=1\ U_{12}=1\ U_{22}=1\ U_{13}=0$ | Infeasible |
| Node 11. $U_{11}=1\ U_{21}=1\ U_{12}=1\ U_{22}=1\ U_{13}=1\ U_{23}=1$ | $\begin{cases} z=257-10\cdot\theta_1, & 0\le\theta_1\le 1.2 \\ z=245, & 1.2\le\theta_1\le 10 \end{cases}$ |

This branch is terminated, and a new upper bound is established, because node 11's optimal solution satisfies all of the requirements of the original mp-MILP model. Because all branches have terminated at this point, node 11's solution is the final optimal solution.

$$\begin{cases} \omega_1=\left[6+10\cdot\theta_1,5,10,12-10\cdot\theta_1,10,10,1,1,1,1,1,1,0,0,0,0\right] \\ z_1=257-10\cdot\theta_1, \quad 0\le\theta_1\le 1.2 \\ \omega_2=\left[18,5,10,0,10,10,1,1,1,1,1,1,0,0,0,0\right] \\ z_2=245, \quad 1.2\le\theta_1\le 10 \end{cases} \tag{B.19}$$

The number of nodes that must be computed using our technique is depicted in Figure 5. This problem has ten binary variables; using the traditional B&B method, the nodes of $2^0+2^1+\cdots+2^{10}$ will be determined. This study's algorithm requires only 11 nodes to be calculated, considerably cutting the number of nodes to be calculated. Figure 5 depicts the exploring map for all potential nodes based on our approach.

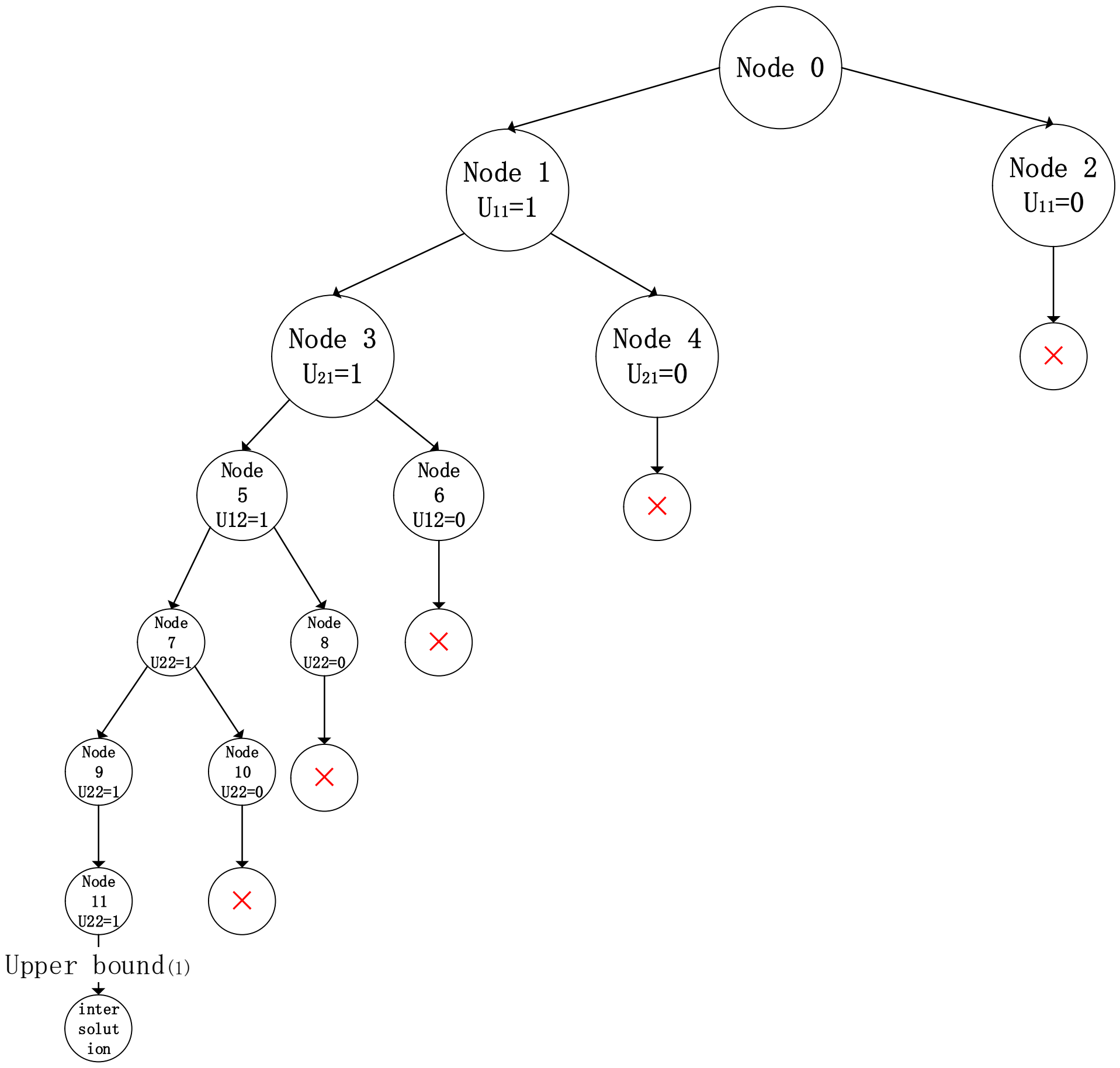

**Figure 5: Nodes based on B&B method**